# HÖRMANDER TYPE THEOREM FOR MULTILINEAR PSEUDO-DIFFERENTIAL OPERATORS

YARYONG HEO, SUNGGEUM HONG AND CHAN WOO YANG

ABSTRACT. We establish a Hörmander type theorem for the multilinear pseudo-differential operators, which is also a generalization of the results in [32] to symbols depending on the spatial variable. Most known results for multilinear pseudo-differential operators were obtained by assuming their symbols satisfy pointwise derivative estimates(Mihlin-type condition), that is, their symbols belong to some symbol classes $n$-$\mathscr{S}^m_{\rho,\delta}(\mathbb{R}^d)$, $0 \le \delta \le \rho \le 1$, $0 \le \delta < 1$ for some $m \le 0$. In this paper, we shall consider multilinear pseudo-differential operators whose symbols have limited smoothness described in terms of function space and not in a pointwise form(Hörmander type condition). Our conditions for symbols are weaker than the Mihlin-type conditions in two senses: the one is that we only assume the first-order derivative conditions in the spatial variable and lower-order derivative conditions in the frequency variable, and the other is that we make use of $L^2$-average condition rather than pointwise derivative conditions for the symbols. As an application, we obtain some mapping properties for the multilinear pseudo-differential operators associated with symbols belonging to the classes $n$-$\mathscr{S}^m_{\rho,\delta}(\mathbb{R}^d)$, $0 \le \rho \le 1$, $0 \le \delta < 1$, $m \le 0$. Moreover, it can be pointed out that our results can be applied to wider classes of symbols which do not belong to the traditional symbol classes $n$-$\mathscr{S}^m_{\rho,\delta}(\mathbb{R}^d)$.

## 1. INTRODUCTION AND STATEMENT OF MAIN RESULTS

Let $\mathscr{S}(\mathbb{R}^d)$ denote the collection of Schwartz functions on $\mathbb{R}^d$, and let $n$ be a positive integer greater than 1. We associate a bounded function $\mathbf{m}(\cdot,\vec{\cdot})$ on $\mathbb{R}^d \times (\mathbb{R}^d)^n$ with $n$-linear pseudo-differential operator $\mathrm{T}_{\mathbf{m}}$ defined by

$$\mathrm{T}_{\mathbf{m}}\big(f_1,\cdots,f_n\big)(x) := \int_{(\mathbb{R}^d)^n} e^{2\pi i\langle x,\xi_1+\cdots+\xi_n\rangle}\,\mathbf{m}(x,\vec{\xi})\widehat{f_1}(\xi_1)\cdots\widehat{f_n}(\xi_n)\,d\vec{\xi}$$

where $f_1,\cdots,f_n$ are Schwartz functions on $\mathbb{R}^d$, $\vec{\xi} := (\xi_1,\dots,\xi_n) \in (\mathbb{R}^d)^n$, and $d\vec{\xi} := d\xi_1\cdots d\xi_n$. Here, $\widehat{f}$ denotes the Fourier transform of $f \in \mathscr{S}(\mathbb{R}^d)$.

Mapping properties of these operators when the symbol $\mathbf{m}$ is independent of the spatial variable $x$, have been well understood in various articles([4, 15–20, 32, 39, 45]). Especially, when $\mathbf{m}$ is independent of the spatial variable $x$, Lee et al. [32] obtained almost sharp results for $H^{p_1}\times\cdots\times H^{p_n} \to L^p$ boundedness of $\mathrm{T}_{\mathbf{m}}$ under the Hörmander type multiplier condition

$$\sup_{j\in\mathbb{Z}} \left\|\mathbf{m}(2^j\,\vec{\cdot}\,)\widehat{\Psi}(\,\vec{\cdot}\,)\right\|_{L^2_s((\mathbb{R}^d)^n)} < \infty, \tag{1.1}$$

where $H^p$ denotes the Hardy space, $\widehat{\Psi}$ is a smooth function that is supported in $\{1/2 < |\vec{\xi}| < 2\}$, and $L^2_s$ denotes $L^2$- Sobolev norm of order $s$.

*Date*: April 14, 2023.

2010 *Mathematics Subject Classification.* Primary 42B20, Secondary 42B15.

This research was supported by Basic Science Research Program through the National Research Foundation of Korea(NRF) funded by the Ministry of Education, Science and Technology NRF-2021R1F1A1045384, NRF-2022R1A4A1018904, and NRF-2022R1A2C1004616, by a Korea University Grant, and by Chosun University(2021).

When the symbol $\mathbf{m}$ depends on the spatial variable $x$, most results for $T_{\mathbf{m}}$ were obtained by assuming the symbol $\mathbf{m}$ satisfy the following Mihlin-type derivative conditions:

$$\left|\partial_x^\alpha \partial_{\vec{\xi}}^\beta \mathbf{m}(x,\vec{\xi})\right| \le C_{\alpha,\beta}(1+|\vec{\xi}|)^{m+\delta|\alpha|-\rho|\beta|} \tag{1.2}$$

for some $0 \le \delta \le \rho \le 1$, $0 \le \delta < 1$, $m \le 0$, and for all multi-indices $\alpha$ and $\beta$. In this case we say that $\mathbf{m}$ belongs to the symbol class $n\text{-}\mathcal{S}^m_{\rho,\delta}(\mathbb{R}^d)$. Unlike the Mihlin-type derivative condition (1.2)(pointwise control), the Hörmander type condition (1.1) contains $L^2$-average of the symbol and its derivatives. Many kinds of symbols satisfy the condition (1.1) rather than (1.2). In this paper, we establish a Hörmander type theorem for the multilinear pseudo-differential operators, which is also a generalization of our previous results [32] to symbols depending on the spatial variable. Under a Hörmander type condition of the symbol $\mathbf{m}$, we establish $h^{p_1}(\mathbb{R}^d)\times\cdots\times h^{p_n}(\mathbb{R}^d) \to L^p(\mathbb{R}^d)$ boundedness of $T_{\mathbf{m}}$ for $0 < p_1, \cdots, p_n < \infty$ and $0 < p < \infty$ satisfying $1/p = 1/p_1 + \cdots + 1/p_n$, where $h^p$ denotes the local Hardy space introduced by Goldberg [11, 12]. Unlike the Mihlin-type symbol condition (1.2), the Hörmander type symbol condition (1.1) can treat various kinds of symbols. Also, while most results for Mihlin-type symbols are obtained by assuming higher-order derivative conditions of the symbols for the spatial variable $x$, our results assume only the first-order derivative conditions of the symbols for the spatial variable $x$. As an application, we obtain some mapping properties for the multilinear pseudo-differential operators associated with the symbols belonging to the classes $n\text{-}\mathcal{S}^m_{\rho,\delta}(\mathbb{R}^d)$, $0 \le \rho \le 1$, $0 \le \delta < 1$, $m \le 0$, as in Theorem 2.1 below. Moreover, it can be pointed out that our results are applied to wider classes of symbols that do not belong to the traditional symbol classes $n\text{-}\mathcal{S}^m_{\rho,\delta}(\mathbb{R}^d)$.

**Remark 1.** In the case of the symbol class $\mathcal{S}^0_{1,0}(\mathbb{R}^d)$, there were some efforts to lower the regularity of the symbol $\mathbf{m}$ in the $x$-variable. Around 1969, L. Nirenberg asked the following question: if we assume that the symbol $\mathbf{m}$ satisfies the derivative condition

$$\left|\partial_\xi^\beta \mathbf{m}(x,\xi)\right| \le C_\beta(1+|\xi|)^{-|\beta|}, \quad x,\xi \in \mathbb{R}^d,$$

for all multi-indices $\beta$, with no a priori regularity in the $x$-variable, does it follow that $T_{\mathbf{m}} : L^2 \to L^2$? In 1972, this question was answered in the negative by Ching [6]. Ching's counterexample is smooth in the $x$-variable, but its $x$-derivatives lack a pointwise control as in (1.2). Meanwhile, in 1978, Coifman-Meyer [8] proved that $T_{\mathbf{m}} : L^p \to L^p$ for $1 < p < \infty$, when $\mathbf{m}$ satisfies the following mild regularity conditions in the $x$-variable

$$\begin{gathered} \left|\partial_\xi^\beta \mathbf{m}(x,\xi)\right| \le C_\beta(1+|\xi|)^{-|\beta|}, \\ \left|\partial_\xi^\beta \mathbf{m}(x+h,\xi) - \partial_\xi^\beta \mathbf{m}(x,\xi)\right| \le C_\beta\,\omega(|h|)(1+|\xi|)^{-|\beta|}, \quad \int_0^1 \omega^2(t)\frac{dt}{t} < \infty. \end{gathered} \tag{1.3}$$

For the bilinear case, under these type conditions (1.3), Coifman-Meyer [8], Maldonado-Naibo [35], and Maldonado [34] obtained $L^p \times L^q \to L^r$ boundedness of $T_{\mathbf{m}}$, for all $1 < p, q < \infty$, $1/p+1/q = 1/r$.

Before we state the main theorem, we first present some known results for multilinear operator $T_{\mathbf{m}}$. To do this we divide the results into two cases: the one is the case where the symbol $\mathbf{m}$ is independent of the space variable $x$ and the other is the complementary case. We contain the results for the first case and the second case in Subsections 1.1 and 1.2, respectively.

1.1. **The case where the symbol m is independent of $x$:** In this case, we abuse the notation to write $\mathbf{m}(x,\vec{\xi}) = \mathbf{m}(\vec{\xi})$. Since we often use the fractional Sobolev spaces to describe previous results and to state our main theorem, we precisely define them here.

For $s \geq 0$ let $(\vec{I}-\vec{\Delta})^{s/2}$ denote the inhomogeneous fractional Laplacian operator acting on functions on $(\mathbb{R}^d)^n$. To be specific,

$$(\vec{I}-\vec{\Delta})^{s/2}F = \Big(\big(1+4\pi^2(|\cdot_1|^2+\cdots+|\cdot_n|^2)\big)^{s/2}\widehat{F}\Big)^{\vee}$$

for a function $F$ on $(\mathbb{R}^d)^n$, where $f^{\vee}(\xi) := \widehat{f}(-\xi)$ denotes the inverse Fourier transform. Now for $s \geq 0$ and $0 < r < \infty$ we define the Sobolev norm

$$\|F\|_{L^r_s((\mathbb{R}^d)^n)} := \big\|(\vec{I}-\vec{\Delta})^{s/2}F\big\|_{L^r((\mathbb{R}^d)^n)}.$$

For the special case $r = 2$, it can be written in the form

$$\|F\|_{L^2_s((\mathbb{R}^d)^n)} = \Big(\int_{(\mathbb{R}^d)^n} \big(1+4\pi^2(|\xi_1|^2+\cdots+|\xi_n|^2)\big)^s \big|\widehat{F}(\xi_1,\dots,\xi_n)\big|^2 d\vec{\xi}\Big)^{1/2}.$$

We first take account of the case $n = 1$, that is, when the operator $T_{\mathbf{m}}$ is a linear operator associated with a multiplier $\mathbf{m}(\xi)$. In this case, the operator $T_{\mathbf{m}}$ in the above can be written as

$$T_{\mathbf{m}}f(x) := \int_{\mathbb{R}^d} e^{2\pi i\langle x,\xi\rangle}\mathbf{m}(\xi)\widehat{f}(\xi)d\xi$$

for $f \in \mathscr{S}(\mathbb{R}^d)$. By Plancherel's identity, we first have $\|T_{\mathbf{m}}\|_{L^2(\mathbb{R}^d)\to L^2(\mathbb{R}^d)} = \|\mathbf{m}\|_{L^\infty(\mathbb{R}^d)}$. According to the classical Mihlin multiplier theorem [36], the operator $T_{\mathbf{m}}$ admits the $L^p$-bounded extension for $1 < p < \infty$ whenever

$$\big|\partial_\xi^\alpha \mathbf{m}(\xi)\big| \leq C_\alpha|\xi|^{-|\alpha|}, \quad \xi \neq 0 \tag{1.4}$$

for all multi-indices $\alpha$ with $|\alpha| \leq [d/2]+1$, and this result was refined by Hörmander [26] who replaced (1.4) with the weaker condition

$$\sup_{j\in\mathbb{Z}}\big\|\mathbf{m}(2^j\cdot)\widehat{\psi}(\cdot)\big\|_{L^2_s(\mathbb{R}^d)} < \infty \quad \text{for } s > d/2, \tag{1.5}$$

where $L^2_s(\mathbb{R}^d)$ stands for the fractional Sobolev space on $\mathbb{R}^d$ and $\psi$ is a Schwartz function on $\mathbb{R}^d$ whose Fourier transform $\widehat{\psi}$ is supported in the annulus $\{\xi \in \mathbb{R}^d : 1/2 < |\xi| < 2\}$ and satisfies $\sum_{j\in\mathbb{Z}}\widehat{\psi}(2^{-j}\xi) = 1$ for all $\xi \neq 0$.

Calderón-Torchinsky [4] extended this result to the (real) Hardy space $H^p(\mathbb{R}^d)$. More precisely they assumed the same condition as in (1.5) with $s > d/p - d/2$ to obtain that for $0 < p \leq 1$ there exists $C > 0$ such that

$$\big\|T_{\mathbf{m}}\big\|_{H^p(\mathbb{R}^d)\to H^p(\mathbb{R}^d)} \leq C\sup_{j\in\mathbb{Z}}\big\|\mathbf{m}(2^j\cdot)\widehat{\psi}(\cdot)\big\|_{L^2_s(\mathbb{R}^d)}. \tag{1.6}$$

The Hardy space $H^p(\mathbb{R}^d)$ is naturally extended over $p > 1$ so that it coincides with $L^p(\mathbb{R}^d)$ for $1 < p \leq \infty$. Recently, the estimates in (1.6) have been reformulated by Grafakos-He-Honzík-Nguyen [16] in this context, namely, if $s > d/r$ and $s > |d/p - d/2|$ with $1 < p < \infty$ and $1 < r < \infty$, then there exists $C > 0$ such that

$$\big\|T_{\mathbf{m}}\big\|_{L^p(\mathbb{R}^d)\to L^p(\mathbb{R}^d)} \leq C\sup_{j\in\mathbb{Z}}\big\|\mathbf{m}(2^j\cdot)\widehat{\psi}\big\|_{L^r_s(\mathbb{R}^d)}. \tag{1.7}$$

We remark that it can be proven that two conditions $s > d/r$ and $s > |d/p - d/2|$ in the above are sharp in the sense that if one of them does not hold, then there exists a bounded function $\mathbf{m}$ for which (1.7) does not hold(see [22, 43]).

Now we turn our attention to the cases $n \geq 2$, that is, the cases where the operators $T_{\mathbf{m}}$ are multilinear operators associated with the multiplier $\mathbf{m}$. For a bounded function $\mathbf{m}$ on $(\mathbb{R}^d)^n$, the operators $T_{\mathbf{m}}$ in the above are called $n$-linear Fourier multipliers which can be rewritten as

$$T_{\mathbf{m}}(f_1, \cdots, f_n)(x) := \int_{(\mathbb{R}^d)^n} e^{2\pi i \langle x, \xi_1 + \cdots + \xi_n \rangle} \mathbf{m}(\vec{\xi}) \widehat{f_1}(\xi_1) \cdots \widehat{f_n}(\xi_n) \, d\vec{\xi}.$$

As a multilinear extension of Mihlin's result, Coifman-Meyer [7, 8] proved that if $L$ is sufficiently large and $\mathbf{m}$ satisfies

$$\left| \partial_{\xi_1}^{\alpha_1} \cdots \partial_{\xi_n}^{\alpha_n} \mathbf{m}(\xi_1, \ldots, \xi_n) \right| \lesssim_{\alpha_1, \ldots, \alpha_n} \left( |\xi_1| + \cdots + |\xi_n| \right)^{-(|\alpha_1| + \cdots + |\alpha_n|)} \tag{1.8}$$

for multi-indices $\alpha_1, \ldots, \alpha_n$ with $|\alpha_1| + \cdots + |\alpha_n| \leq L$, then $T_{\mathbf{m}}$ is bounded from $L^{p_1}(\mathbb{R}^d) \times \cdots \times L^{p_n}(\mathbb{R}^d)$ to $L^p(\mathbb{R}^d)$ for all $1 < p_1, \ldots, p_n \leq \infty$ and $1 < p < \infty$ with $1/p_1 + \cdots + 1/p_n = 1/p$. The result was extended to the case $p \leq 1$ by Kenig-Stein [30] and Grafakos-Torres [25]. Later, the research naturally proceeded toward improving the condition (1.8) to obtain multilinear analogs of the classical Hörmander multiplier theorem, which was initiated by Tomita in [45], where he considered the $n$-linear counterpart $\Psi$ of $\psi$ in the multilinear context, that is, $\Psi$ is a Schwartz function on $(\mathbb{R}^d)^n$ having the properties that

$$\operatorname{supp}(\widehat{\Psi}) \subset \left\{ \vec{\xi} \in (\mathbb{R}^d)^n : 1/2 \leq |\vec{\xi}| \leq 2 \right\}, \qquad \sum_{j \in \mathbb{Z}} \widehat{\Psi}(2^{-j} \vec{\xi}) = 1, \quad \vec{\xi} \neq \vec{0}$$

and obtained that if for $1 < p, p_1, \ldots, p_n < \infty$, $1/p = 1/p_1 + \cdots + 1/p_n$ and

$$\sup_{j \in \mathbb{Z}} \left\| \mathbf{m}(2^j \vec{\cdot}) \widehat{\Psi}(\vec{\cdot}) \right\|_{L^2_s((\mathbb{R}^d)^n)} < \infty \tag{1.9}$$

with $s > nd/2$, then

$$\left\| T_{\mathbf{m}} \right\|_{L^{p_1} \times \cdots \times L^{p_n} \to L^p} \lesssim \sup_{j \in \mathbb{Z}} \left\| \mathbf{m}(2^j \vec{\cdot}) \widehat{\Psi}(\vec{\cdot}) \right\|_{L^2_s((\mathbb{R}^d)^n)}. \tag{1.10}$$

This was extended by Grafakos-Si [23] to the range $p \leq 1$ in terms of the $L^r$-based Sobolev space condition for $1 < r \leq 2$. Later, the standard Sobolev spaces in the estimates (1.10) have been replaced by product-type Sobolev spaces in many recent results. For $s_1, \ldots, s_n \geq 0$, we define the product-type Sobolev spaces $L^2_{(s_1, \ldots, s_n)}((\mathbb{R}^d)^n)$ as function spaces consisting of all functions $F$ on $(\mathbb{R}^d)^n$ such that the norm

$$\|F\|_{L^2_{(s_1, \ldots, s_n)}((\mathbb{R}^d)^n)} := \left( \int_{(\mathbb{R}^d)^n} \left(1 + 4\pi^2 |\xi_1|^2\right)^{s_1} \cdots \left(1 + 4\pi^2 |\xi_n|^2\right)^{s_n} \left| \widehat{F}(\xi_1, \ldots, \xi_n) \right|^2 d\vec{\xi} \right)^{1/2}$$

is finite. Miyachi-Tomita in [39] replaced the condition (1.9) with the condition of the product-type Sobolev spaces

$$\sup_{j \in \mathbb{Z}} \left\| \mathbf{m}(2^j \vec{\cdot}) \widehat{\Psi}(\vec{\cdot}) \right\|_{L^2_{(s_1, \ldots, s_n)}((\mathbb{R}^d)^n)} < \infty,$$

to obtain $H^{p_1}(\mathbb{R}^d) \times H^{p_2}(\mathbb{R}^d) \to L^p(\mathbb{R}^d)$ boundedness for bilinear multipliers (i.e., $n = 2$) in the full range of indices $0 < p, p_1, p_2 \leq \infty$ extending the estimates in (1.6) to the bilinear setting. Multilinear extensions were later provided by Grafakos-Miyachi-Tomita [18], Grafakos-Nguyen [20], and Grafakos-Miyachi-Nguyen-Tomita [19]. One can combine results in [17–20, 39] to present them in one formulation as follows:
Let $0 < p_1, \ldots, p_n \leq \infty$, $0 < p < \infty$, and $1/p_1 + \cdots + 1/p_n = 1/p$. Suppose that

$$s_1, \ldots, s_n > \frac{d}{2}, \qquad \sum_{k \in I} \left( \frac{s_k}{d} - \frac{1}{p_k} \right) > -\frac{1}{2} \tag{1.11}$$

for any nonempty subsets $I$ of $J_n := \{1,\dots,n\}$. Then every $\mathrm{T_m}$ satisfies

$$\big\|\mathrm{T}_{\mathbf{m}}\big(f_1,\dots,f_n\big)\big\|_{L^p(\mathbb{R}^d)} \lesssim \sup_{j\in\mathbb{Z}} \big\|\mathbf{m}(2^j\vec{\cdot}\,)\widehat{\Psi}(\vec{\cdot}\,)\big\|_{L^2_{(s_1,\dots,s_n)}((\mathbb{R}^d)^n)} \prod_{i=1}^n \|f_i\|_{H^{p_i}(\mathbb{R}^d)} \tag{1.12}$$

for Schwartz functions $f_1,\dots,f_n$ on $\mathbb{R}^d$.

Now we come back to the original condition (1.9). The necessary conditions in this setting were obtained in [15]. Precisely, it was established that for $0<p,p_1,\dots,p_n<\infty$ with $1/p = 1/p_1+\cdots+1/p_n$ and $0<r,s<\infty$ if we suppose that

$$\|\mathrm{T}_{\mathbf{m}}\|_{L^{p_1}\times\cdots\times L^{p_n}\to L^p} \lesssim \sup_{j\in\mathbb{Z}} \big\|\mathbf{m}(2^j\vec{\cdot}\,)\widehat{\Psi}(\vec{\cdot}\,)\big\|_{L^r_s((\mathbb{R}^d)^n)} \tag{1.13}$$

for all bounded functions $\mathbf{m}$ for which $\sup_{j\in\mathbb{Z}} \big\|\mathbf{m}(2^j\vec{\cdot}\,)\widehat{\Psi}(\vec{\cdot}\,)\big\|_{L^r_s((\mathbb{R}^d)^n)}<\infty$, then it is necessary to have

(1) $s\ge \max\big\{\frac{(n-1)d}{2}, \frac{nd}{r}\big\}$,
(2) $\frac{1}{p}-\frac{1}{2}\le \frac{s}{d}+\sum_{i\in I}\big(\frac{1}{p_i}-\frac{1}{2}\big)$ where $I$ is an arbitrary subset of $J_n=\{1,2,\dots,n\}$ which may also be empty (in which case the sum is supposed to be zero).

Recently, Lee et al. [32] consider the case $r=2$ in (1.13), and they proved that the necessary conditions (1) and (2) in the above are also "almost" sufficient for the $H^{p_1}\times\cdots\times H^{p_n}\to L^p$ boundedness for $\mathrm{T_m}$.

**Theorem A** ([32])**.** Let $\mathbf{m}=\mathbf{m}(\vec{\xi})$. Let $0<p_1,\cdots,p_n\le\infty$ and $0<p<\infty$ satisfy $1/p=1/p_1+\cdots+1/p_n$. Suppose that

(1) $s>\frac{nd}{2}$,
(2) $\frac{1}{p}-\frac{1}{2}< \frac{s}{d}+\sum_{i\in I}\big(\frac{1}{p_i}-\frac{1}{2}\big)$ where $I$ is an arbitrary subset of $J_n=\{1,2,\dots,n\}$ which may also be empty (in which case the sum is supposed to be zero).

Then we have

$$\big\|\mathrm{T}_{\mathbf{m}}(f_1,\dots,f_n)\big\|_{L^p(\mathbb{R}^d)} \le C\sup_{j\in\mathbb{Z}} \big\|\mathbf{m}(2^j\vec{\cdot}\,)\widehat{\Psi}(\vec{\cdot}\,)\big\|_{L^2_s((\mathbb{R}^d)^n)} \prod_{i=1}^n \|f_i\|_{H^{p_i}(\mathbb{R}^d)}, \tag{1.14}$$

for Schwartz functions $f_1,\dots,f_n$ on $\mathbb{R}^d$.

Because of the necessary condition in the above, the conditions (1) and (2) in Theorem A are “almost" sharp except for the critical case

$$s=\frac{nd}{2} \qquad \text{or} \qquad \frac{1}{p}-\frac{1}{2}=\frac{s}{d}+\sum_{i\in I}\Big(\frac{1}{p_i}-\frac{1}{2}\Big) \quad \text{for some } \ I\subset J_n.$$

Also two conditions $s>nd/2$ and $1/p-1/2<s/d$ are necessary for (1.14) to hold. We conclude this subsection by stating a lemma which is an equivalent classification of the condition (2) in Theorem A.

**Lemma 1.1.** *The set of all collection of* $\big(\frac{1}{p_1},\dots,\frac{1}{p_n}\big)\in(0,\infty)^n$ *that satisfies the condition* (2) *in Theorem A is equivalent to the set* $\mathrm{B}_n\big(\frac{s}{d}+\frac{1}{2}\big)$ *where*

$$\mathrm{B}_n(\alpha) := \Big\{(x_1,\cdots,x_n)\in(0,\infty)^n \, : \, \sum_{i=1}^n \max\big(x_i,1/2\big)<\alpha\Big\}. \tag{1.15}$$

*Proof.* The proof will be given in Section 9(Appendix: Proof of Lemma 1.1). □

Now we turn our attention to multilinear multiplier theory for pseudo-differential operators.

1.2. **The case where the symbol m depends on** $x$**:** Compared to the previous case, properties of the multilinear operators $T_{\mathbf{m}}$ associated with symbols depending on the spatial variable $x$ have not been well understood. Most results for $T_{\mathbf{m}}$ were obtained by assuming $\mathbf{m}$ belongs to some symbol classes $n\text{-}\mathcal{S}^m_{\rho,\delta}(\mathbb{R}^d)$, $0 \le \delta \le \rho \le 1$, $0 \le \delta < 1$ for some $m \le 0$. That is,

$$\left|\partial_x^\alpha \partial_{\vec{\xi}}^\beta \mathbf{m}(x,\vec{\xi})\right| \le C_{\alpha,\beta}(1+|\vec{\xi}|)^{m+\delta|\alpha|-\rho|\beta|} \tag{1.16}$$

for all multi-indices $\alpha$ and $\beta$. The number $m$ is called the order of $\mathbf{m}$. For related results for multilinear operator $T_{\mathbf{m}}$, we refer to the following papers: Bényi et al. [1], Bényi-Maldonado-Naibo-Torres [2], Bényi-Torres [3], Coifman-Meyer [7, 8], Huang-Chen [28], Kato-Miyachi-Tomita [29], Nirenberg [31], Lu-Zhang [33], and Miyachi-Tomita [38, 40].

Next, we recall the definition of the local Hardy space $h^p$ introduced by Goldberg [11, 12].

**The real Hardy space** $H^p(\mathbb{R}^d)$**.** First, we recall the definition of the real Hardy space $H^p(\mathbb{R}^d)$ based on Stein's book [44, Chapter III, §1]. A tempered distribution $f$ is in $H^p$ if and only if $\sup_{t>0}|\varphi_t * f| \in L^p$, here $\varphi_t(x) = t^{-d}\varphi(x/t)$, $\varphi \in \mathcal{S}(\mathbb{R}^d)$, $\int \varphi \ne 0$. For each $0 < p < \infty$, there exists an $N > 0$ so that if $B = \{\varphi \in \mathcal{S} : \|\varphi\|_{\alpha,\beta} := \sup_{x\in\mathbb{R}^d}|x^\alpha \partial_x^\beta \varphi(x)| \le 1 \text{ for } |\alpha|, |\beta| \le N\}$, and if $\psi \in \mathcal{S}$ with $\int \psi \ne 0$, then the $L^p$ norms of the following functions are equivalent:

$$\sup_{t>0}|\psi_t * f(x)|, \quad \sup_{t>0}\sup_{\varphi\in B}|\varphi_t * f(x)|, \quad \sup_{t>0}\sup_{|y|<t}|\psi_t * f(x-y)|.$$

Any one of these can be taken to be the $H^p$ norm of $f$, and written by $\|f\|_{H^p}$ for $0 < p < \infty$. As was pointed out by Goldberg [11, 12]: if $0 < p \le 1$, then

1. $H^p$ does not contain $\mathcal{S}$;
2. $H^p$ is not well defined on manifolds;
3. pseudo-differential operators are not bounded on $H^p$(compare this to the result (1.6) of Calderón-Torchinsky [4]).

Because of these problems, Goldberg [11, 12] introduced a space $h^p(\mathbb{R}^d)$ which satisfies (1), (2), and (3) in the positive sense. As for (3), Goldberg obtained the following result:

**Theorem B** ([11, 12])**.** Suppose $\mathbf{m} = \mathbf{m}(x,\xi)$ belongs to $\mathcal{S}^0_{1,0}(\mathbb{R}^d)$, then for $0 < p < \infty$

$$\|T_{\mathbf{m}}(f)\|_{h^p} \le C\|f\|_{h^p}.$$

**The local Hardy space** $h^p(\mathbb{R}^d)$**.** For $0 < p < \infty$, let $h^p(\mathbb{R}^d)$ denote the local Hardy space of Goldberg [11, 12]. That is, a tempered distribution $f$ is in $h^p$ if and only if $\sup_{0<t<1}|\varphi_t * f| \in L^p$, here $\varphi_t(x) = t^{-d}\varphi(x/t)$, $\varphi \in \mathcal{S}(\mathbb{R}^d)$, $\int \varphi \ne 0$. Note that $H^p \subset h^p$. As in [11, 12], for each $0 < p < \infty$, there exists an $N > 0$ so that if $B = \{\varphi \in \mathcal{S} : \|\varphi\|_{\alpha,\beta} := \sup_{x\in\mathbb{R}^d}|x^\alpha \partial_x^\beta \varphi(x)| \le 1 \text{ for } |\alpha|, |\beta| \le N\}$, and if $\psi \in \mathcal{S}$ with $\int \psi \ne 0$, then the $L^p$ norms of the following functions are equivalent:

$$\sup_{0<t<1}|\psi_t * f(x)|, \quad \sup_{0<t<1}\sup_{\varphi\in B}|\varphi_t * f(x)|, \quad \sup_{0<t<1/2}\sup_{|y|<t}|\psi_t * f(x-y)|.$$

Any one of these can be taken to be the $h^p$ norm of $f$, and written by $\|f\|_{h^p}$ for $0 < p < \infty$. Note that $H^p = h^p = L^p$ if $1 < p < \infty$, and if $0 < p \le 1$, $f \in L^1_{loc}$, then

$$\|f\|_{L^p} \lesssim \|f\|_{h^p} \lesssim \|f\|_{H^p}.$$

**Some known results for linear case.** We first take account of the case $n = 1$, that is, when the operator $T_{\mathbf{m}}$ is a linear operator associated with a symbol $\mathbf{m}(x,\xi)$. Hörmander [27] proved that if $\mathbf{m} \in \mathscr{S}^0_{\rho,\delta}(\mathbb{R}^d)$, $0 \le \delta < \rho \le 1$, then $T_{\mathbf{m}}$ is bounded on $L^2(\mathbb{R}^d)$. After that, Calderón-Vaillancourt [5] gerenalized this result to the case $0 \le \delta = \rho < 1$. On the other hand, Fefferman [10] showed that if $\mathbf{m} \in \mathscr{S}^{-\mu}_{\rho,\delta}(\mathbb{R}^d)$, $0 \le \delta < \rho \le 1$, $\mu \ge (1-\rho)|\frac{1}{p} - \frac{1}{2}|d$, then $T_{\mathbf{m}}$ is bounded on $L^p(\mathbb{R}^d)$, $1 < p < \infty$. After that, Päivärinta-Somersalo [41] extended this results up to $0 \le \rho = \delta < 1$ and $0 < p \le 1$ with the local Hardy pace $h^p$ introduded by Goldberg [11, 12]. For the special case $\rho = \delta = 0$, these results imply the following.

**Theorem C** ([5, 8, 10, 29, 37, 41])**.** Let $0 < p_1 \le p \le \infty$. Suppose $\mathbf{m} = \mathbf{m}(x,\xi)$ belongs to $\mathscr{S}^m_{0,0}(\mathbb{R}^d)$ for some $m \le 0$. Then

$$\|T_{\mathbf{m}}(f)\|_{h^p} \le C\|f\|_{h^{p_1}},$$

holds if and only if

$$m \le \min\left\{\frac{d}{p}, \frac{d}{2}\right\} - \max\left\{\frac{d}{p_1}, \frac{d}{2}\right\},$$

where $h^p$ ($h^{p_1}, resp.$) should be replaced by $bmo$ when $p = \infty$ ($p_1 = \infty, resp.$).

In case $p = p_1$, the result for Theorem C is given by [5, 8, 10, 37, 41], and as was pointed out in Kato-Miyachi-Tomita [29], the case $p_1 < p$ can be deduced from the case $p = p_1$ with the aid of the mapping properties of the fractional integration operator and symbolic calculus in $\mathscr{S}^m_{0,0}(\mathbb{R}^d)$.

**Some known results for multilinear case.** For the bilinear case, Bényi-Maldonado-Naibo-Torres [2] proved that if $\mathbf{m}$ belongs to the symbol class 2-$\mathscr{S}^0_{1,\delta}(\mathbb{R}^d)$, $0 \le \delta < 1$, then $T_{\mathbf{m}}$ has a bounded extension from $L^p \times L^q$ into $L^r$, for all $1 < p, q < \infty$, $1/p + 1/q = 1/r$.

**Remark 2.** By Bényi et al. [1], if $0 \le \rho < 1$, $0 \le \delta \le 1$, and $1 \le p, q, r < \infty$ such that $1/p + 1/q = 1/r$, then there exist symbols in 2-$\mathscr{S}^0_{\rho,\delta}(\mathbb{R}^d)$ that give rise to unbounded operators from $L^p(\mathbb{R}^d) \times L^q(\mathbb{R}^d)$ into $L^r(\mathbb{R}^d)$.

On the other hand, when $\mathbf{m}$ belongs to the symbol class 2-$\mathscr{S}^m_{\rho,\delta}(\mathbb{R}^d)$, $0 \le \delta \le \rho \le 1$, $\delta < 1$ for some $m < 0$, Bényi et al. [1] obtained the following results.

**Theorem D** ([1])**.** If $\mathbf{m}$ belongs to the symbol class 2-$\mathscr{S}^m_{\rho,\delta}(\mathbb{R}^d)$, $0 \le \delta \le \rho \le 1$, $\delta < 1$, then $T_{\mathbf{m}}$ is bounded from $L^p \times L^q$ to $L^r$ $(1 \le p, q \le \infty, 1/p + 1/q = 1/r)$ under the condition

$$m < m(p,q) := d(\rho - 1)\Big(\max\left\{\frac{1}{2}, \frac{1}{p}, \frac{1}{q}, 1 - \frac{1}{r}\right\} + \max\left\{\frac{1}{r} - 1, 0\right\}\Big).$$

Multilinear cases were considered by Coifman-Meyer [7, 8], they proved that if $\mathbf{m}$ belongs to the symbol class $n$-$\mathscr{S}^0_{1,0}(\mathbb{R}^d)$, then $T_{\mathbf{m}}$ has a bounded extension from $L^{p_1} \times \cdots \times L^{p_n}$ into $L^p$, for all $1 < p_1, \ldots, p_n < \infty$ and $1 \le p < \infty$ with $1/p_1 + \cdots + 1/p_n = 1/p$. After that Grafakos-Torres [25] and Kenig-Stein [30] extended this result up to the optimal range of $p > 1/n$.

Recently, Kato-Miyachi-Tomita [29] obtained the following mapping properties for the multilinear pseudo-differential operators associated with the symbols belonging to the classes $n$-$\mathscr{S}^m_{0,0}(\mathbb{R}^d)$ for $m \le 0$. (For $n = 2$, the results were obtained by Miyachi-Tomita [38], and generalized to the case $n \ge 2$ by Kato-Miyachi-Tomita [29].)

**Theorem E** ([29, 38])**.** Let $\mathbf{m} = \mathbf{m}(x,\vec{\xi})$ belong to $n$-$\mathscr{S}^m_{0,0}(\mathbb{R}^d)$ for some $m \in \mathbb{R}$. Let $n \ge 2$, $0 < p, p_1, \cdots, p_n \le \infty$ and $1/p \le 1/p_1 + \cdots + 1/p_n$. Then

$$\big\|T_{\mathbf{m}}(f_1, \ldots, f_n)\|_{h^p(\mathbb{R}^d)} \le C \prod_{i=1}^n \|f_i\|_{h^{p_i}(\mathbb{R}^d)}, \tag{1.17}$$

holds if and only if

$$m \le \min\left\{\frac{d}{p}, \frac{d}{2}\right\} - \sum_{i=1}^{n} \max\left\{\frac{d}{p_i}, \frac{d}{2}\right\}. \tag{1.18}$$

If (1.18) is satisfied and if some of $p_i$'s, $1 \le i \le n$, are equal to $\infty$, then (1.17) holds with the corresponding $h^{p_i}$ replaced by $bmo$.

To state our main theorem, we let $\Phi$ be a Schwartz function on $(\mathbb{R}^d)^n$ whose Fourier transform $\widehat{\Phi}$ is supported in $|\vec{\xi}| < 1$ and $\widehat{\Phi}(\vec{\xi}) = 1$ for $|\vec{\xi}| \le 1/2$. Together with $\Phi$, we define another function $\Psi$ by $\widehat{\Psi}(\vec{\xi}) = \widehat{\Phi}(\vec{\xi}) - \widehat{\Phi}(2\vec{\xi})$. Then we have the following "partition of unity" of the $\vec{\xi}$-space:

$$1 = \widehat{\Phi}(\vec{\xi}) + \sum_{j=0}^{\infty} \widehat{\Psi}(2^{-j}\vec{\xi}), \quad \text{for all } \vec{\xi}.$$

Note that $\widehat{\Psi}$ is supported in the annulus of the form $\{\vec{\xi} : 1/2 < |\vec{\xi}| < 2\}$.

We now state the main theorem:

**Theorem 1.2.** *Let $n \ge 1$. Let $\mathrm{B}_n(\alpha)$ be as in (1.15). Let $\mathbf{m} = \mathbf{m}(x, \vec{\xi})$ be a bounded function on $\mathbb{R}^d \times (\mathbb{R}^d)^n$. Let $0 < p_1, \cdots, p_n < \infty$ and $0 < p < \infty$ satisfy $1/p = 1/p_1 + \cdots + 1/p_n$. Suppose that*

*(1) $s > \frac{nd}{2}$,*

*(2) $\mathrm{B}_n(\frac{s}{d}) = \left\{(x_1, \cdots, x_n) \in (0, \infty)^n : \sum_{i=1}^n \max(x_i, 1/2) < \frac{s}{d}\right\}$.*

*Then for any $0 \le \delta < 1$, if $\left(\frac{1}{p_1}, \ldots, \frac{1}{p_n}\right) \in \mathrm{B}_n(\frac{s}{d})$, then*

$$\left\|\mathrm{T}_{\mathbf{m}}(f_1, \ldots, f_n)\right\|_{L^p(\mathbb{R}^d)} \le C_{s,\delta} \, \|\mathbf{m}\|_{\mathscr{L}^2_{s,\delta}} \prod_{i=1}^{n} \|f_i\|_{h^{p_i}(\mathbb{R}^d)}, \tag{1.19}$$

*for Schwartz functions $f_1, \ldots, f_n$ on $\mathbb{R}^d$, where*

$$\begin{aligned}\|\mathbf{m}\|_{\mathscr{L}^2_{s,\delta}} := \sup_{x \in \mathbb{R}^d} &\left(\sum_{|\alpha| \le 1} \left\|\partial_x^\alpha \mathbf{m}(x, \vec{\cdot}\,)\widehat{\Phi}(\vec{\cdot}\,)\right\|_{L^2_s((\mathbb{R}^d)^n)}\right) \\ &+ \sup_{j \ge 0} \sup_{x \in \mathbb{R}^d} \left(\sum_{|\alpha| \le 1} 2^{-j\delta|\alpha|} \left\|\partial_x^\alpha \mathbf{m}(x, 2^j \vec{\cdot}\,)\widehat{\Psi}(\vec{\cdot}\,)\right\|_{L^2_s((\mathbb{R}^d)^n)}\right).\end{aligned} \tag{1.20}$$

We illustrate the domains $\mathrm{B}_n(\frac{s}{d})$ in Figures 1 and 2 for $n = 2$ and 3, respectively.

**Remark 3.**
(1) The condition $\|\mathbf{m}\|_{\mathscr{L}^2_{s,\delta}} < \infty$ in (1.20) is a natural generalization of the symbol class $n$-$\mathscr{S}^m_{\rho,\delta}(\mathbb{R}^d)$ in (1.16). Moreover, we only assume at most the first-order differentiability of the symbol concerning the space variable $x \in \mathbb{R}^d$ in (1.20).

(2) In the previous paper [32], we considered the case where $\mathbf{m}$ does not depend on $x$ and obtained the strong type estimates

$$\left\|\mathrm{T}_{\mathbf{m}}(f_1, \ldots, f_n)\right\|_{L^p(\mathbb{R}^d)} \le C_s \left(\sup_{j \in \mathbb{Z}} \left\|\mathbf{m}(2^j \vec{\cdot}\,)\widehat{\Psi}(\vec{\cdot}\,)\right\|_{L^2_s((\mathbb{R}^d)^n)}\right) \prod_{i=1}^{n} \|f_i\|_{H^{p_i}(\mathbb{R}^d)},$$

for Schwartz functions $f_1, \ldots, f_n$ on $\mathbb{R}^d$ if $\left(\frac{1}{p_1}, \ldots, \frac{1}{p_n}\right) \in \mathrm{B}_n(\frac{s}{d} + \frac{1}{2})$ and $s > \frac{nd}{2}$.

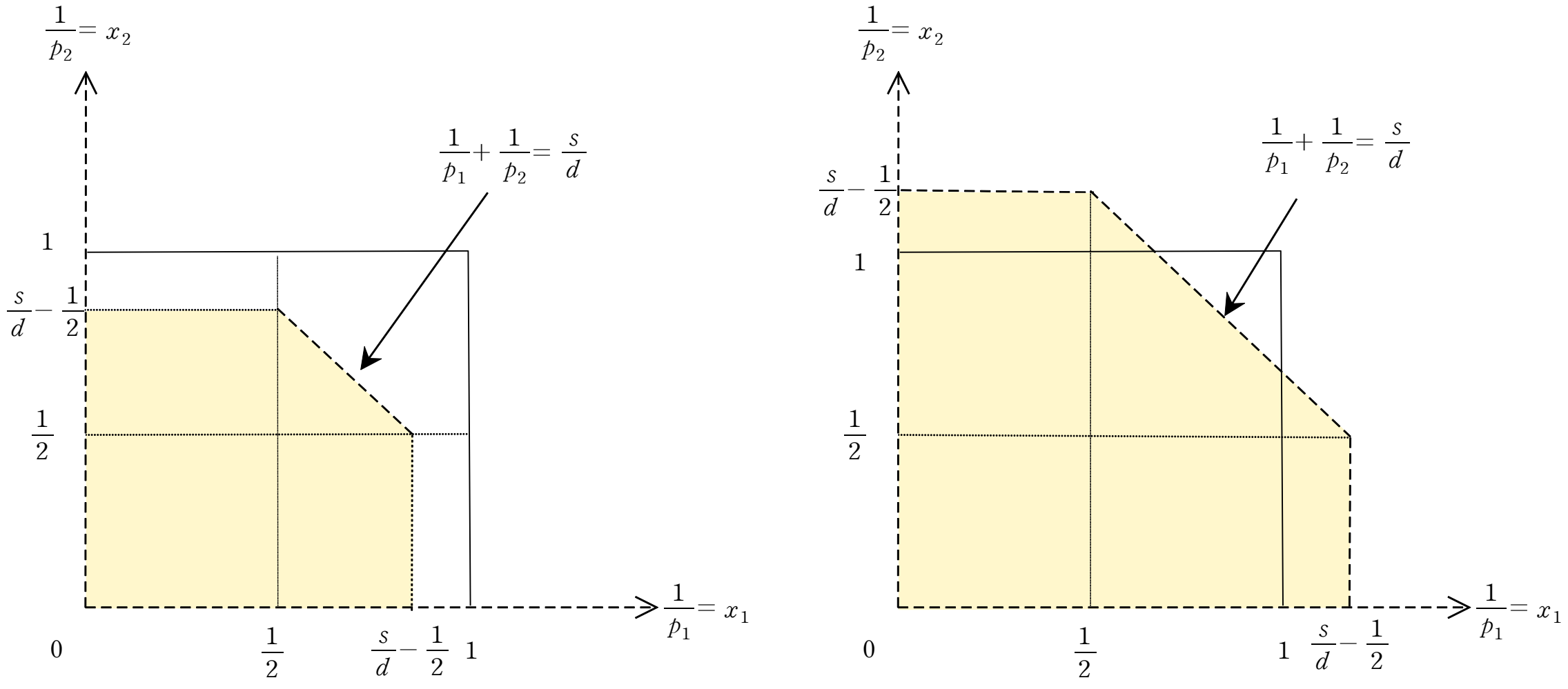

(a) $B_2(\frac{s}{d})$ for $1 < \frac{s}{d} \le \frac{3}{2}$ (b) $B_2(\frac{s}{d})$ for $\frac{3}{2} < \frac{s}{d}$

FIGURE 1. $B_n(\frac{s}{d})$ for $n = 2$

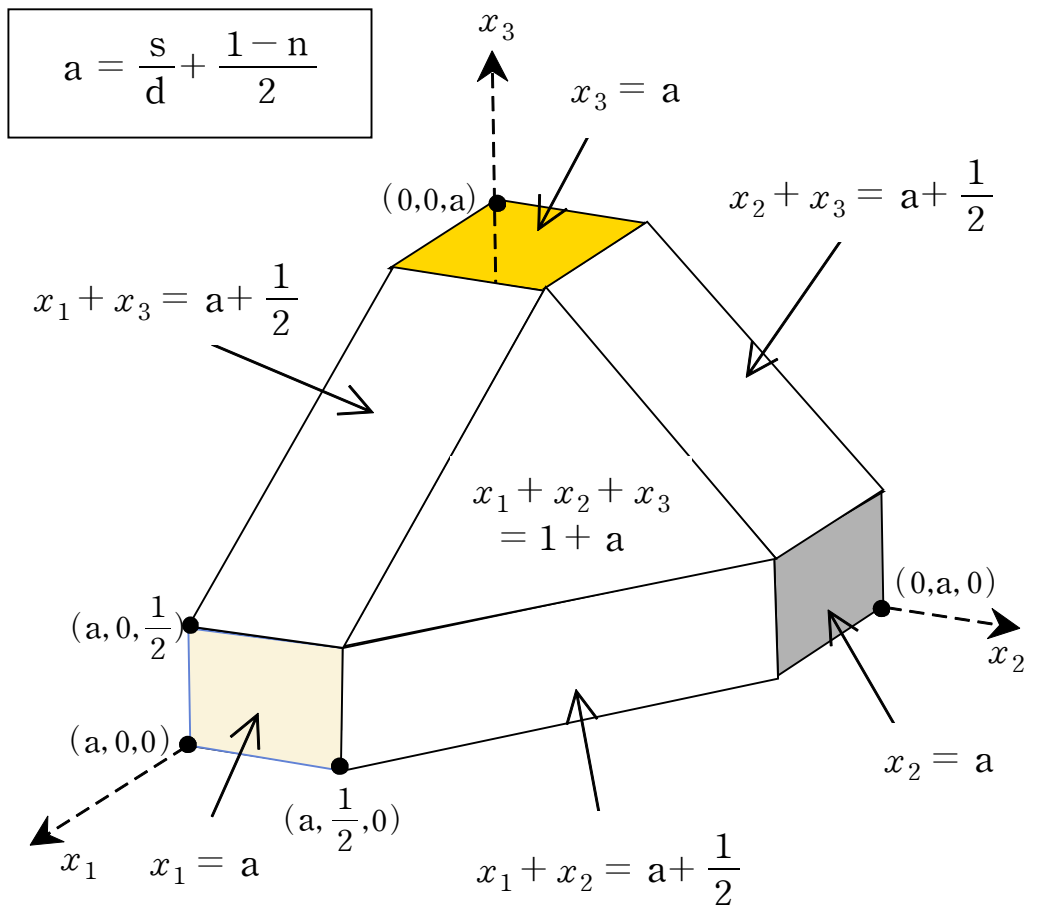

FIGURE 2. $B_n(\frac{s}{d})$ for $n = 3$

(3) Although Theorem A in [32] holds for $\left(\frac{1}{p_1}, \ldots, \frac{1}{p_n}\right) \in B_n(\frac{s}{d} + \frac{1}{2})$, we obtain Theorem 1.2 for $\left(\frac{1}{p_1}, \ldots, \frac{1}{p_n}\right) \in B_n(\frac{s}{d})$. At present we do not know whether or not our results can be extended to $B_n(\frac{s}{d} + \frac{1}{2})$. Let us briefly explain the reason. In the previous paper [32] we used something like duality arguments $\big(\langle Tf, g\rangle = \langle f, T^*g\rangle\big)$: that is, note that

$$\langle T_{\mathbf{m}}(f_1, \cdots, f_n), f_{n+1}\rangle = \int_{\mathbb{R}^d} \int_{(\mathbb{R}^d)^n} \mathscr{F}^{-1}[\mathbf{m}(\vec{\cdot})](\vec{y}) \Big( \prod_{i=1}^{n} f_i(x - y_i) \Big) f_{n+1}(x) d\vec{y} dx.$$

Then by using translation $x \to x + y_1$ and Hölder's inequality we have

$$\begin{aligned}&|\langle \mathrm{T}_{\mathbf{m}}(f_1, \cdots, f_n), f_{n+1}\rangle| \\ &\le \left(\|\mathscr{F}^{-1}[\mathbf{m}(\vec{\cdot})](\vec{y})\|_{L^2(d\vec{y})}\right)\int_{\mathbb{R}^d}|f_1(x)|\left\|\Big(\prod_{i=2}^n f_i(x-(y_i-y_1))\Big)f_{n+1}(x+y_1)\right\|_{L^2(d\vec{y})}dx.\end{aligned} \tag{1.21}$$

The important thing in (1.21) is the following: by using translation $x \to x + y_1$ we control the right hand side of (1.21) by the $L^1$ integral of $f_1$ (instead of $L^2$), which allows us to gain $\frac{1}{2}$ on $\mathrm{B}_n(\frac{s}{d}+\frac{1}{2})$. But in the case $\mathbf{m}$ depends on $x$, if we use the translation $x \to x + y_1$, the term $\|\mathscr{F}^{-1}[\mathbf{m}(\vec{\cdot})](\vec{y})\|_{L^2(d\vec{y})}$ in (1.21) becomes $\sup_x \left\|\mathscr{F}^{-1}[\mathbf{m}(x+y_1, \vec{\cdot})](\vec{y})\right\|_{L^2(d\vec{y})}$, and our arguments used in the previous paper [32] break down at this point.

## 2. Examples and applications

In this section, we give various examples that satisfy the Hörmander type symbol condition (1.20), and explain advantages of using Hörmander type symbol condition.

Recall that, by Bényi et al. [1], if $0 \le \rho < 1$, $0 \le \delta \le 1$, and $1 \le p, q, r < \infty$ such that $1/p + 1/q = 1/r$, then there exist symbols in $2$-$\mathscr{S}^0_{\rho,\delta}(\mathbb{R}^d)$ that give rise to unbounded operators from $L^p(\mathbb{R}^d) \times L^q(\mathbb{R}^d)$ into $L^r(\mathbb{R}^d)$. So it is natural to consider the symbol classes $n$-$\mathscr{S}^m_{\rho,\delta}(\mathbb{R}^d)$, $(0 \le \rho \le 1, 0 \le \delta < 1)$ for some $m \le 0$. In the case $0 \le \delta \le \rho \le 1$, $\delta < 1$, and $m < 0$, there is a positive results for bilinear operators(see Theorem D).

Our main result in Theorem 1.2 can be applied to the symbol classes $n$-$\mathscr{S}^m_{\rho,\delta}(\mathbb{R}^d)$, $(0 \le \rho \le 1, 0 \le \delta < 1, m \le 0)$ to obtain the following results for multilinear operators. Note that the condition in Theorem 2.1 below does not depend on the order of $\delta$ and $\rho$.

**Theorem 2.1.** *Let $0 < p_1, \cdots, p_n < \infty$ and $0 < p < \infty$ satisfy $1/p = 1/p_1 + \cdots + 1/p_n$. Let $0 \le \rho \le 1$, $0 \le \delta < 1$, and $m \le 0$. Suppose that $\mathbf{m}(x, \vec{\xi})$ satisfies the derivative conditions*

$$\left|\partial_x^\alpha \partial_{\vec{\xi}}^\beta \mathbf{m}(x, \vec{\xi})\right| \le C_{\alpha,\beta}(1+|\vec{\xi}|)^{m+\delta|\alpha|-\rho|\beta|} \tag{2.1}$$

*for all multi-indices $\alpha$ and $\beta$ with $|\alpha| \le 1$. If $m \le (\rho-1)s$, $s > \frac{nd}{2}$, and $\left(\frac{1}{p_1}, \ldots, \frac{1}{p_n}\right) \in \mathrm{B}_n(\frac{s}{d})$, then*

$$\left\|\mathrm{T}_{\mathbf{m}}(f_1, \ldots, f_n)\right\|_{L^p(\mathbb{R}^d)} \le C \prod_{i=1}^n \|f_i\|_{h^{p_i}(\mathbb{R}^d)}, \tag{2.2}$$

*for Schwartz functions $f_1, \ldots, f_n$ on $\mathbb{R}^d$.*

*Proof of Theorem 2.1.* When $s$ is a non-negative integer, then by using the derivative conditions (2.1) we have

$$\sup_x \left\|\partial_x^\alpha \mathbf{m}(x, \vec{\cdot})\widehat{\Phi}(\vec{\cdot})\right\|_{L^2_s((\mathbb{R}^d)^n)} \lesssim 1, \quad \sup_x \left\|\partial_x^\alpha \mathbf{m}(x, 2^j\vec{\cdot})\widehat{\Psi}(\vec{\cdot})\right\|_{L^2_s((\mathbb{R}^d)^n)} \lesssim 2^{j\delta|\alpha|}2^{j(1-\rho)s+jm}. \tag{2.3}$$

Let $s \ge 0$ be a real number, then choose non-negative integer $\nu$ such that $\nu \le s \le \nu+1$. Then by interpolating the results of two cases $\nu$ and $\nu+1$, we see that (2.3) holds for $s$. By (2.3), if $m + (1-\rho)s \le 0$, then $\|\mathbf{m}\|_{\mathscr{L}^2_{s,\delta}} < \infty$. And the result follows from Theorem 1.2. □

Theorem 2.1 implies the following corollary.

**Corollary 2.2.** *Let $0 \le \rho \le 1$, $0 \le \delta < 1$. If $\mathbf{m}$ belongs to the symbol classes $n$-$\mathscr{S}^m_{\rho,\delta}(\mathbb{R}^d)$ with $m \le (\rho-1)s$ for some $s > \frac{nd}{2}$. Then (2.2) holds if $\left(\frac{1}{p_1}, \ldots, \frac{1}{p_n}\right) \in \mathrm{B}_n(\frac{s}{d})$. In particular, we have the following.*

*(1) If* $\mathbf{m}$ *belongs to the symbol classes* $n\text{-}\mathscr{S}^0_{1,\delta}(\mathbb{R}^d)$*, then by taking* $s \to \infty$*, we see that* (2.2) *holds for any* $0 < p, p_1, \cdots, p_n < \infty$ *satisfying* $1/p = 1/p_1 + \cdots + 1/p_n$.

*(2) Suppose that* $0 \le \rho < 1$, $\frac{-m}{(1-\rho)} > \frac{nd}{2}$*. If* $\mathbf{m}$ *belongs to the symbol classes* $n\text{-}\mathscr{S}^m_{\rho,\delta}(\mathbb{R}^d)$*, then* (2.2) *holds if* $\left(\frac{1}{p_1}, \ldots, \frac{1}{p_n}\right) \in \mathrm{B}_n(\frac{-m}{(1-\rho)d})$.

**Remark 4.**
(1) The condition $\left(\frac{1}{p_1}, \ldots, \frac{1}{p_n}\right) \in \mathrm{B}_n(\frac{-m}{(1-\rho)d})$ in Corollary 2.2 (2) can be written as

$$\sum_{i=1}^{n} \max\left\{\frac{1}{p_i}, \frac{1}{2}\right\} < \frac{-m}{(1-\rho)d}. \tag{2.4}$$

By Corollary 2.2, our main Theorem 1.2, when it is applied to the symbol classes $n\text{-}\mathscr{S}^m_{\rho,\delta}(\mathbb{R}^d)$, $(0 \le \rho \le 1, 0 \le \delta < 1)$, implies better results as $\rho \nearrow 1$.

(2) Although most results for $\mathrm{T}_{\mathbf{m}}$ were obtained by assuming $\mathbf{m}$ belongs to some symbol classes $n\text{-}\mathscr{S}^m_{\rho,\delta}(\mathbb{R}^d)$, $0 \le \delta \le \rho \le 1, 0 \le \delta < 1$ for some $m \le 0$, our results contain the case $\rho < \delta$ if $m < \frac{(\rho-1)nd}{2}$. To the best of our knowledge, we are the first to obtain the results on this case.

(3) Note that, the condition (1.18) in Theorem E is equivalent to

$$\sum_{i=1}^{n} \max\left\{\frac{1}{p_i}, \frac{1}{2}\right\} \le \frac{-m}{d} + \min\left\{\frac{1}{p}, \frac{1}{2}\right\}. \tag{2.5}$$

By (2.4) and (2.5), Theorem 1.2, when applied to the symbol classes $n\text{-}\mathscr{S}^m_{0,0}(\mathbb{R}^d)$, implies a weaker result than that of Kato et al. [29](Theorem E). But our theorem can be applied to any symbol classes $n\text{-}\mathscr{S}^m_{\rho,\delta}(\mathbb{R}^d)$, $(0 \le \rho \le 1, 0 \le \delta < 1)$, and use only the first-order derivative conditions of the symbol for the spatial variable $x$. As far as we know, except for the symbol class $n\text{-}\mathscr{S}^0_{1,0}(\mathbb{R}^d)$, including [29], most studies on multilinear pseudo-differential operators associated with the symbol classes $n\text{-}\mathscr{S}^m_{\rho,\delta}(\mathbb{R}^d)$, $0 \le \delta \le \rho < 1$, $m \le 0$, have been carried out under higher order derivative conditions of the symbol concerning the spatial variable $x$.

Now we take some examples of symbols that do not belong to the traditional symbols $n\text{-}\mathscr{S}^m_{\rho,\delta}(\mathbb{R}^d)$.

**Example 1.** Let $\phi(x, \vec{\xi})$ be a smooth function that is supported in $\left\{(x, \vec{\xi}) \in \mathbb{R}^d \times (\mathbb{R}^d)^n : |(x, \vec{\xi})| < 1\right\}$. Let $\varphi$ be a smooth function. We consider the multiplier

$$\mathbf{m}(x, \vec{\xi}) = \phi(x, \vec{\xi}) e^{i|x|^{\frac{3}{2}}\varphi(\vec{\xi})}.$$

Then

$$\left|\partial_x^\alpha \partial_{\vec{\xi}}^\beta \mathbf{m}(x, \vec{\xi})\right| \le C_{\alpha,\beta}$$

for all multi-indices $\beta$, and $\alpha$ only for $|\alpha| \le 1$. Since $\|\mathbf{m}\|_{\mathscr{L}^2_{s,0}} < \infty$ for any $s > 0$, (2.2) holds for all $0 < p, p_1, \cdots, p_n < \infty$ satisfying $1/p = 1/p_1 + \cdots + 1/p_n$.

**Example 2** (cf. [24])**.** Unlike the Mihlin-type symbol condition, the Hörmander type condition can treat symbols whose derivatives have infinitely many singularities. Let $0 \le \delta < 1$. For positive integers $k$, let $\Phi_k(x, \vec{\xi})$ be a function defined on $\mathbb{R}^d \times (\mathbb{R}^d)^n$. Suppose that $\Phi_k(x, \vec{\xi})$ are supported in $\{\vec{\xi} \in (\mathbb{R}^d)^n : \frac{1}{2} < |\vec{\xi}| < 2\}$ and satisfy the derivative conditions

$$\left|\partial_x^\alpha \partial_{\vec{\xi}}^\beta \Phi_k(x, \vec{\xi})\right| \le C_{\alpha,\beta} 2^{k\delta|\alpha|}$$

for all multi-indices $\alpha$ and $\beta$ with $|\alpha| \le 1$. Let $\vec{a}_k \in \{\vec{\xi} \in (\mathbb{R}^d)^n : \frac{1}{2} < |\vec{\xi}| < 2\}$. We consider the multiplier

$$\mathbf{m}(x,\vec{\xi}) = \sum_{k=1}^{\infty} \Phi_k(x, 2^{-k}\vec{\xi})\, |2^{-k}\vec{\xi} - \vec{a}_k|^{\gamma} \tag{2.6}$$

for some $\gamma > 0$. Then

$$\|\mathbf{m}\|_{\mathscr{L}^2_{s,\delta}} < \infty \quad \text{if } s < \gamma + \frac{nd}{2}.$$

To see this, by (2.6) it suffices to show that

$$\sup_x \left\| \partial_x^{\alpha} \Phi_k(x,\vec{\xi}) |\vec{\xi} - \vec{a}_k|^{\gamma} \right\|_{L^2_s} \le C_{\alpha} 2^{k\delta|\alpha|} \quad (|\alpha| \le 1) \tag{2.7}$$

if $0 \le s < \gamma + \frac{nd}{2}$. Let $\phi$ be a smooth function that is supported in $\{1/2 < |\vec{\xi}| < 2\}$ and

$$\sum_{l=-4}^{\infty} \phi(2^l \vec{\xi}) = 1 \quad \text{if } 0 < |\vec{\xi}| < 4.$$

Then by (2.7) it suffices to show that

$$\sup_x \left\| \partial_x^{\alpha} \Phi_k(x,\vec{\xi}) \phi(2^l(\vec{\xi} - \vec{a}_k)) |\vec{\xi} - \vec{a}_k|^{\gamma} \right\|_{L^2_s} \le C_{\alpha} 2^{k\delta|\alpha|} \Big( 2^{l(-\gamma + s - \frac{nd}{2})} \Big). \tag{2.8}$$

(2.8) is clear when $s$ is a non-negative integer. Let $s \ge 0$ be a real number, then choose non-negative integer $\nu$ such that $\nu \le s \le \nu + 1$. The result for $s$ follows by interpolating the results of two cases $\nu$ and $\nu+1$. Therefore if $\big(\frac{1}{p_1}, \dots, \frac{1}{p_n}\big) \in \mathrm{B}_n(\frac{\gamma}{d} + \frac{n}{2})$, then (2.2) holds. Note that $\mathbf{m}$ does not belong to any symbol class $n\text{-}\mathscr{S}^m_{\rho,\delta}(\mathbb{R}^d)$.

**Example 3.** Let $0 \le \delta < 1$. For positive integers $k$, let $\phi_k(x,\vec{\xi})$ be a function defined on $\mathbb{R}^d \times (\mathbb{R}^d)^n$. Suppose that $\phi_k(x,\vec{\xi})$ are supported in $\{\vec{\xi} \in (\mathbb{R}^d)^n : \frac{1}{2} < |\vec{\xi}| < 2\}$ and satisfy the derivative conditions

$$\left| \partial_x^{\alpha} \partial_{\vec{\xi}}^{\beta} \phi_k(x,\vec{\xi}) \right| \le C_{\alpha,\beta} 2^{k\delta|\alpha|}$$

for all multi-indices $\alpha$ and $\beta$ with $|\alpha| \le 1$. Let $\vec{a}_k(x) : \mathbb{R}^d \to (\mathbb{R}^d)^n$ such that

$$\left| \partial_x^{\alpha} \vec{a}_k(x) \right| \le C_{\alpha} 2^{k\delta|\alpha|}$$

for all multi-indices $\alpha$ with $|\alpha| \le 1$. We consider the multiplier

$$\mathbf{m}(x,\vec{\xi}) = \sum_{k=1}^{\infty} \phi_k(x, 2^{-k}\vec{\xi})\, |2^{-k}\vec{\xi} - \vec{a}_k(x)|^{\gamma}$$

for some $\gamma > 1$. As in Example 2 we have $\|\mathbf{m}\|_{\mathscr{L}^2_{s,\delta}} < \infty$ if $s < (\gamma - 1) + \frac{nd}{2}$. Therefore if $\big(\frac{1}{p_1}, \dots, \frac{1}{p_n}\big) \in \mathrm{B}_n(\frac{\gamma-1}{d} + \frac{n}{2})$, then (2.2) holds.

**Example 4** (cf. [10]). Let $0 \le \delta < 1$. For positive integers $k$, let $\phi_k(x,\vec{\xi})$ be a function defined on $\mathbb{R}^d \times (\mathbb{R}^d)^n$. Suppose that $\phi_k(x,\vec{\xi})$ are supported in $\{\vec{\xi} \in (\mathbb{R}^d)^n : \frac{1}{2} < |\vec{\xi}| < 2\}$ and satisfy the derivative conditions

$$\left| \partial_x^{\alpha} \partial_{\vec{\xi}}^{\beta} \phi_k(x,\vec{\xi}) \right| \le C_{\alpha,\beta} 2^{k\delta|\alpha|} \tag{2.9}$$

for all multi-indices $\alpha$ and $\beta$ with $|\alpha| \le 1$. For positive real numbers $a$ and $b$, we consider the multiplier

$$\mathbf{m}(x,\vec{\xi}) = \sum_{k=1}^{\infty} \phi_k(x, 2^{-k}\vec{\xi}) |\vec{\xi}|^{-b} e^{i|\vec{\xi}|^a}.$$

If $0 \le a \le 1$, and (2.9) holds for all $\alpha$, then $\mathbf{m} \in n\text{-}\mathscr{S}^{-b}_{1-a,\delta}(\mathbb{R}^d)$. But if $a > 1$, then $\mathbf{m}$ does not belong to any symbol class $n\text{-}\mathscr{S}^{m}_{\rho,\delta}(\mathbb{R}^d)$. Note that $\|\mathbf{m}\|_{\mathscr{L}^2_{s,\delta}} < \infty$ if $-b + as \le 0$. Therefore if $\frac{nd}{2} < s \le \frac{b}{a}$ and $\big(\frac{1}{p_1}, \dots, \frac{1}{p_n}\big) \in \mathrm{B}_n(\frac{s}{d})$, then (2.2) holds. This means that if $\frac{b}{a} > \frac{nd}{2}$ and $\big(\frac{1}{p_1}, \dots, \frac{1}{p_n}\big) \in \mathrm{B}_n(\frac{b}{ad})$, then (2.2) holds.

## 3. Preliminaries

In the remaining part of this paper, we make use of the following notations.

**Notation 1.** *We use the notation $\langle \cdot, \cdot \rangle$ to denote both the inner product of functions and the dot product of points. That is, $\langle f, g\rangle = \int_{\mathbb{R}^d} f(x)\overline{g(x)}\,dx$ for two functions $f$ and $g$, and $\langle a, b\rangle = a \cdot b$ for two points $a, b \in \mathbb{R}^d$. For two quantities $A$ and $B$, we shall write $A \lesssim B$ if $A \le CB$ holds for some positive constant $C$, depending on the dimension and possibly other parameters apparent from the context. We write $A \sim B$ if both $A \lesssim B$ and $B \lesssim A$ hold. For a measurable set $E$, the notation $|E|$ stands for the measure of $E$ and $\chi_E$ does the characteristic function of $E$. The symbol $\sharp S$ means the cardinality of the set $S$.*

In this section, we establish several estimates which will be used in the rest sections of this paper.

**Lemma 3.1.** *For $0 < r < \infty$, let $\mathrm{M}_r f := \big(\mathrm{M}(|f|^r)\big)^{1/r}$ where $\mathrm{M}$ denotes the Hardy-Littlewood maximal operator. Then, using the Fefferman-Stein vector-valued maximal inequality in [9], we obtain that*

$$\Big\|\big\{\mathrm{M}_r(f_j)\big\}_{j\in\mathbb{Z}}\Big\|_{L^p(\ell^q)} \lesssim \Big\|\big\{f_j\big\}_{j\in\mathbb{Z}}\Big\|_{L^p(\ell^q)},$$

*provided that $0 < p < \infty$, $0 < q \le \infty$, and $0 < r < p, q$.*

**Lemma 3.2** (Grafakos-Si [23])**.** *Let $\triangle_j$ be the Littlewood-Paley operator given by $\widehat{\triangle_j(g)}(\xi) = \widehat{g}(\xi)\widehat{\Psi}(2^{-j}\xi)$, $j \in \mathbb{Z}$. Suppose that a tempered distribution $f$ satisfies*

$$\Big\|\Big(\sum_{j\in\mathbb{Z}} |\triangle_j(f)|^2\Big)^{\frac{1}{2}}\Big\|_{L^p} < \infty,$$

*and the support of $\widehat{f} \subset \mathbb{R}^d \setminus \{0\}$. Then for $0 < p < \infty$*

$$\|f\|_{H^p(\mathbb{R}^d)} \le c(d, p, \Psi)\Big\|\Big(\sum_{j\in\mathbb{Z}} |\triangle_j(f)|^2\Big)^{1/2}\Big\|_{L^p(\mathbb{R}^d)}$$

*where $H^p(\mathbb{R}^d)$ denotes the Hardy space on $\mathbb{R}^d$. For the proof see Lemma 2.4 in Grafakos-Si [23], or Theorem 2.2.9 in [14].*

**Lemma 3.3** ([14], Theorem 2.2.9)**.** *Let $\Psi \in \mathscr{S}(\mathbb{R}^d)$ whose Fourier transform is compactly supported away from the origin. For each $j \in \mathbb{Z}$, let $\Psi_j(\cdot) := 2^{jd}\Psi(2^j\cdot)$. Let $0 < p < \infty$. Then for all $f \in H^p(\mathbb{R}^d)$ we have*

$$\Big\|\Big(\sum_{j\in\mathbb{Z}} |f * \Psi_j|^2\Big)^{\frac{1}{2}}\Big\|_{L^p(\mathbb{R}^d)} \le C_{d,p,\Psi}\|f\|_{H^p(\mathbb{R}^d)}.$$

**Lemma 3.4** ([11, 14])**.** *Let $\Psi \in \mathscr{S}(\mathbb{R}^d)$ whose Fourier transform is compactly supported away from the origin. For each $j \in \mathbb{Z}$, let $\Psi_j(\cdot) := 2^{jd}\Psi(2^j\cdot)$. Let $0 < p < \infty$. Then for all $f \in h^p(\mathbb{R}^d)$ we have*

$$\Big\|\Big(\sum_{j=0}^{\infty} |f * \Psi_j|^2\Big)^{\frac{1}{2}}\Big\|_{L^p(\mathbb{R}^d)} \le C_{d,p,\Psi}\|f\|_{h^p(\mathbb{R}^d)}.$$

*Proof.* The proof will be given in Section 9(Appendix: Proof of Lemma 3.4). For the proof, we adopt the proof of Theorem 2.2.9 in [14] and [11, Theorem B]. □

The following Lemma is taken from Grafakos's book [14, Lemma 2.2.3].

**Lemma 3.5** ([14], Lemma 2.2.3). *Let $0 < r < \infty$. Then there exists a constant $C$ such that for all $t > 0$ and for all $\mathscr{C}^1$ functions $g$ on $\mathbb{R}^d$ whose distributional Fourier transform is supported in the ball $|\xi| \le t$ we have*

$$\sup_{z\in\mathbb{R}^d} \frac{|g(x-z)|}{(1+t|z|)^{\frac{d}{r}}} \le C\big(\mathrm{M}(|g|^r)(x)\big)^{\frac{1}{r}}, \tag{3.1}$$

*where* M *denotes the Hardy-Littlewood maximal operator. The constant $C$ depends only on the dimension $d$ and $r$; in particular, it is independent of $t$.*

**Lemma 3.6.** *Let $\Phi \in \mathscr{S}(\mathbb{R}^d)$ whose Fourier transform is compactly supported. Let $M \ge 0$. For any $0 < r < \infty$*

$$\sup_{|y|\lesssim 1} \big|[(f*\Phi_t)](x+t2^M y)\big| \lesssim 2^{\frac{Md}{r}} \mathrm{M}_r[(f*\Phi_t)](x),$$

*where $\mathrm{M}_r f := \big(\mathrm{M}(|f|^r)\big)^{1/r}$ and $\Phi_t(x) = t^{-d}\Phi(x/t)$.*

*Proof of Lemma 3.6.* Since $\widehat{(f*\Phi_t)}(\xi) = \widehat{f}(\xi)\widehat{\Phi}(t\xi)$ is supported in $|\xi| \lesssim t^{-1}$, the result follows by applying (3.1) with $z = t2^M y$. □

**Lemma 3.7.** *Let $\Phi \in \mathscr{S}(\mathbb{R}^d)$ whose Fourier transform is compactly supported and $\Phi_k(x) = 2^{kd}\Phi(2^k x)$. Let $M \ge 0$ and $k \ge j$, then for any $0 < r \le s$ we have*

$$\Big(\int_{|y|\lesssim 1} \big|f*\Phi_k(x+2^{-j+M}y)\big|^s dy\Big)^{\frac{1}{s}} \lesssim 2^{(k-j+M)d(\frac{1}{r}-\frac{1}{s})}\mathrm{M}_r[f*\Phi_k](x). \tag{3.2}$$

*Proof of Lemma 3.7.* Since $\widehat{f*\Phi_k}(\xi) = \widehat{f}(\xi)\widehat{\Phi}(2^{-k}\xi)$ is supported in $|\xi| \lesssim 2^k$, by Lemma 3.6 with $t = 2^{-k}$, for any $0 < r \le s$ we have

$$\begin{aligned} &\int_{|y|\lesssim 1} \big|f*\Phi_k(x+2^{-j+M}y)\big|^s dy \\ &\lesssim \Big(\sup_{|y|\lesssim 1} \big|(f*\Phi_k)(x+2^{-j+M}y)\big|^{s-r}\Big)\Big(\int_{|y|\lesssim 1} \big|(f*\Phi_k)(x+2^{-j+M}y)\big|^r dy\Big) \\ &\lesssim \Big(2^{(k-j+M)\frac{d}{r}}\mathrm{M}_r[(f*\Phi_k)](x)\Big)^{s-r}\Big(\mathrm{M}(|(f*\Phi_k)|^r)(x)\Big) \\ &= 2^{(k-j+M)d(\frac{s}{r}-1)}\Big(\mathrm{M}_r[(f*\Phi_k)](x)\Big)^s. \end{aligned} \tag{3.3}$$

By (3.3) we have (3.2). □

**Lemma 3.8.** *Let $\Phi \in \mathscr{S}(\mathbb{R}^d)$ whose Fourier transform is compactly supported and $\Phi_j(x) = 2^{jd}\Phi(2^j x)$. Let $M \ge 0$. Then for any $0 < r \le s$ and $r < p < \infty$ we have*

$$\Big\|\sup_{j\ge 0}\Big(\int_{|y|\lesssim 1} \big|f*\Phi_j(x+2^{-j+M}y)\big|^s dy\Big)^{\frac{1}{s}}\Big\|_{L^p(\mathbb{R}^d)} \lesssim 2^{Md\left(\frac{1}{r}-\frac{1}{s}\right)}\|f\|_{h^p(\mathbb{R}^d)}. \tag{3.4}$$

*Proof of Lemma 3.8.* By Lemma 3.7, for any $0 < r \le s$ we have

$$\sup_{j\ge 0}\Big(\int_{|y|\lesssim 1} \big|f*\Phi_j(x+2^{-j+M}y)\big|^s dy\Big)^{\frac{1}{s}} \lesssim 2^{Md\left(\frac{1}{r}-\frac{1}{s}\right)}\sup_{j\ge 0}\mathrm{M}_r[f*\Phi_j](x).$$

Then for $r < p < \infty$

$$\|\sup_{j\ge 0}\mathrm{M}_r[f*\Phi_j]\|_{L^p(\mathbb{R}^d)} \le \|\mathrm{M}_r[\sup_{j\ge 0}|f*\Phi_j|]\|_{L^p(\mathbb{R}^d)} \lesssim \|\sup_{j\ge 0}|f*\Phi_j|\|_{L^p(\mathbb{R}^d)} \lesssim \|f\|_{h^p(\mathbb{R}^d)},$$

and we have (3.4). □

**Lemma 3.9.** *Let $\Phi \in \mathscr{S}(\mathbb{R}^d)$ whose Fourier transform is compactly supported and $\Phi_k(x) = 2^{kd}\Phi(2^k x)$. Let $\omega^N(x) = (1+|x|)^{-N}$ and $\omega_k^N(x) = 2^{kd}\omega^N(2^k x)$. Let $M \geq 0$ and $k \geq j$. Let $1 \leq s < \infty$. For any $0 < r \leq s$, if $N > (\frac{s}{r}+1)d$, then*

$$\left(\int_{|y|\lesssim 1}\Big[\omega_k^N * \big|f*\Phi_k\big|(x+2^{-j+M}y)\Big]^s dy\right)^{\frac{1}{s}} \lesssim 2^{(k-j+M)d(\frac{1}{r}-\frac{1}{s})}\mathrm{M}_r\big[(f*\Phi_k)\big](x). \tag{3.5}$$

*Proof of Lemma 3.9.* If $N > d$, then by Hölder's inequality we have

$$\Big[\omega_k^N * \big|f*\Phi_k\big|(x)\Big]^s \lesssim \omega_k^N * |f*\Phi_k|^s(x). \tag{3.6}$$

We write

$$\begin{aligned}
&\omega_k^N * |f*\Phi_k|^s(x+2^{-j+M}y)\\
&= \int_{\mathbb{R}^d} |f*\Phi_k|^s\big(x+2^{-j+M}y-z\big)\frac{2^{kd}}{(1+|2^k z|)^N}dz\\
&= \int_{\mathbb{R}^d} |f*\Phi_k|^s\big(x+2^{-j+M}y-2^{-j+M}z'\big)\frac{2^{(k-j+M)d}}{(1+|2^{k-j+M}z'|)^N}dz'\\
&= \int_{\mathbb{R}^d} \big|(f*\Phi_k)\big|^s\big(x+2^{-k}(2^{k-j+M}(y-z'))\big)\frac{(2^{k-j+M})^d}{(1+2^{k-j+M}|z'|)^N}dz'\\
&\lesssim \int_{|z'|\leq 1}\big(\cdots\big)dy' + \sum_{l\geq 0}\int_{|z'|\sim 2^l}\big(\cdots\big)dy'.
\end{aligned} \tag{3.7}$$

Let $|z'| \sim 2^l$ for $l \geq 0$, then $|y-z'| \lesssim 2^l$ and by change of variable $y - z' \to y$ we have

$$\begin{aligned}
&\int_{|y|\lesssim 1}\int_{|z'|\sim 2^l} \big|(f*\Phi_k)\big|^s\big(x+2^{-k}(2^{k-j+M}(y-z'))\big)\frac{(2^{k-j+M})^d}{(1+2^{k-j+M}|z'|)^N}dz'dy\\
&\lesssim \left(\frac{(2^{k-j+M})^d}{(1+2^{k-j+M+l})^N}2^{ld}\right)\int_{|y|\lesssim 2^l}\big|(f*\Phi_k)\big|^s\big(x+2^{-k}(2^{k-j+M}y)\big)dy.
\end{aligned} \tag{3.8}$$

By Lemma 3.6, for $0 < r \leq s$ and $l \geq 0$ we have

$$\begin{aligned}
&\int_{|y|\lesssim 2^l}\big|(f*\Phi_k)\big|^s\big(x+2^{-k}(2^{k-j+M}y)\big)du\\
&\lesssim \Big[\sup_{|y|\lesssim 2^l}\big|(f*\Phi_k)\big|\big(x+2^{-k}(2^{k-j+M}y)\big)\Big]^{s-r}\int_{|y|\lesssim 2^l}\big|(f*\Phi_k)\big|^r\big(x+2^{-k}(2^{k-j+M}y)\big)dy\\
&\lesssim \Big[(2^{k-j+M+l})^{\frac{d}{r}}\mathrm{M}_r\big[(f*\Phi_k)\big](x)\Big]^{s-r}\Big[2^{ld}\mathrm{M}\big[\big|(f*\Phi_k)\big|^r\big](x)\Big].
\end{aligned} \tag{3.9}$$

By (3.8) and (3.9)

$$\text{the left-hand side of (3.8)} \lesssim (2^{k-j+M})^{\frac{s}{r}d-N}2^{l(\frac{s}{r}d+d-N)}\Big(\mathrm{M}_r\big[(f*\Phi_k)\big](x)\Big)^s. \tag{3.10}$$

Similarly, we have

$$\begin{aligned}
&\int_{|y|\lesssim 1}\int_{|z'|\lesssim 1}\big|(f*\Phi_k)\big|^s\big(x+2^{-k}(2^{k-j+M}(y-z'))\big)\frac{(2^{k-j+M})^d}{(1+2^{k-j+M}|z'|)^N}dz'dy\\
&\lesssim (2^{k-j+M})^{\frac{d}{r}(s-r)}\Big(\mathrm{M}_r(f*\Phi_k)(x)\Big)^s.
\end{aligned} \tag{3.11}$$

By (3.7), (3.10) and (3.11), if $N > (\frac{s}{r}+1)d$, then we have

$$\int_{|y|\lesssim 1} \omega_k^N * |f*\Phi_k|^s(x+2^{-j+M}y)dy \lesssim (2^{k-j+M})^{\frac{d}{r}(s-r)}\Big(\mathrm{M}_r(f*\Phi_k)(x)\Big)^s,$$

and (3.5) follows by applying (3.6) if $1\le s<\infty$. □

**Lemma 3.10.** *Let $0<p,q<\infty$. Suppose $f_N\to f$ in $L^q(\mathbb{R}^d)$ and $\|f_N\|_{L^p(\mathbb{R}^d)}\le A<\infty$ where $A$ is independent of $N$. Then $\|f\|_{L^p(\mathbb{R}^d)}\le 2A$.*

*Proof.* Note that

$$\begin{aligned}\Big|\{x\in\mathbb{R}^d:|f(x)|>\alpha\}\Big| &\le \Big|\{x\in\mathbb{R}^d:|f_N(x)|>\alpha/2\}\Big| + \Big|\{x\in\mathbb{R}^d:|(f-f_N)(x)|>\alpha/2\}\Big|\\ &\le \Big|\{x\in\mathbb{R}^d:|f_N(x)|>\alpha/2\}\Big| + \frac{2^q}{\alpha^q}\|f-f_N\|_{L^q}^q,\end{aligned}$$

and

$$\|f\|_{L^p(\mathbb{R}^d)}^p = p\int_0^\infty \alpha^{p-1}\Big|\{x\in\mathbb{R}^d:|f(x)|>\alpha\}\Big|d\alpha. \tag{3.12}$$

Let $0<\epsilon<M<\infty$, then

$$\begin{aligned} &p\int_\epsilon^M \alpha^{p-1}\Big|\{x\in\mathbb{R}^d:|f(x)|>\alpha\}\Big|d\alpha\\ &\le p\int_0^\infty \alpha^{p-1}\Big|\{x\in\mathbb{R}^d:|f_N(x)|>\alpha/2\}\Big|d\alpha + p\int_\epsilon^M \alpha^{p-1}\Big(\frac{2^q}{\alpha^q}\|f-f_N\|_{L^q}^q\Big)d\alpha\\ &= 2^p\|f_N\|_p^p + p\int_\epsilon^M \alpha^{p-1}\Big(\frac{2^q}{\alpha^q}\|f-f_N\|_{L^q}^q\Big)d\alpha\\ &\le 2^pA^p + p\int_\epsilon^M \alpha^{p-1}\Big(\frac{2^q}{\alpha^q}\|f-f_N\|_{L^q}^q\Big)d\alpha \to 2^pA^p\end{aligned}$$

as $N\to\infty$. Thus we have

$$p\int_\epsilon^M \alpha^{p-1}\Big|\{x\in\mathbb{R}^d:|f(x)|>\alpha\}\Big|d\alpha \le 2^pA^p$$

which is independent of $0<\epsilon<M<\infty$. By (3.12), this implies that $\|f\|_{L^p}\le 2A$. □

**Littlewood-Paley type decomposition of $\mathrm{T}_{\mathbf{m}}$.** Recall that $\psi$ is a Schwartz function on $\mathbb{R}^d$ generating Littlewood-Paley functions $\{\psi_j\}_{j\in\mathbb{Z}}$ with $\mathrm{supp}(\widehat{\psi})\subset\{\xi\in\mathbb{R}^d:1/2\le|\xi|\le 2\}$ and $\sum_{j\in\mathbb{Z}}\widehat{\psi_j}(\xi)=1$ for $\xi\neq 0$ where $\psi_j:=2^{jd}\psi(2^j\cdot)$. Such a function $\psi$ can be constructed as follows. Let $\varphi\in\mathscr{S}(\mathbb{R}^d)$ be a Schwartz function such that

$$\mathrm{supp}(\widehat{\varphi})\subseteq[-2,2]^d \quad\text{and}\quad \widehat{\varphi}(\xi)=1 \text{ on } [-1,1]^d.$$

Then define $\psi\in\mathscr{S}(\mathbb{R}^d)$ so that $\widehat{\psi}(\xi):=\widehat{\varphi}(\xi)-\widehat{\varphi}(2\xi)$. Note that $\mathrm{supp}(\widehat{\psi})\subseteq\{\xi:1/2\le|\xi|\le 2\}$. For each $k\in\mathbb{Z}$ define $\widehat{\psi_k}(\xi):=\widehat{\psi}(2^{-k}\xi)$. Then $\mathrm{supp}(\widehat{\psi_k})\subseteq\{\xi:2^{k-1}\le|\xi|\le 2^{k+1}\}$ and

$$\sum_{k\in\mathbb{Z}}\widehat{\psi_k}(\xi)=1 \quad \text{if } \xi\neq 0. \tag{3.13}$$

The following decomposition lemma is taken from [32]. The essentially same decomposition, which has a different presentation, is described in [21, 42].

**Lemma 3.11** ([32], Lemma 4.1). *Let $\Psi$ be a Schwartz function whose Fourier transform $\widehat{\Psi}$ is supported in $\{\vec{\xi} \in (\mathbb{R}^d)^n : 1/2 \le |\vec{\xi}| \le 2\}$ and satisfies*

$$\widehat{\Phi}(\vec{\xi}) + \sum_{j\ge 0} \widehat{\Psi}(2^{-j}\vec{\xi}) = 1 \quad \text{for all } \vec{\xi} \neq 0.$$

*Then the term $\sum_{j\ge0}\sum_{k_1,k_2,\cdots,k_n\in\mathbb{Z}} \widehat{\Psi}(2^{-j}\vec{\xi})\widehat{\psi_{k_1}}(\xi_1)\widehat{\psi_{k_2}}(\xi_2)\cdots\widehat{\psi_{k_n}}(\xi_n)$ can be written as a finite sum of form*

$$\sum_{j\ge0} \widehat{\Psi}(2^{-j}\vec{\xi})\widehat{\Phi}^1_j(\xi_1)\widehat{\Phi}^2_j(\xi_2)\cdots\widehat{\Phi}^n_j(\xi_n)\widehat{\Phi}^{n+1}_j(-\xi_1-\cdots-\xi_n),$$

*where $\vec{\xi} = (\xi_1, \xi_2, \cdots, \xi_n)$, and $\widehat{\Phi}^1$, $\widehat{\Phi}^2$, $\cdots$, $\widehat{\Phi}^{n+1}$ are compactly supported smooth functions, and at least two of $\widehat{\Phi}^1$, $\widehat{\Phi}^2$, $\cdots$, $\widehat{\Phi}^{n+1}$ are compactly supported away from the origin, and $\widehat{\Phi}^i_j(\cdot) := \widehat{\Phi}^i(2^{-j}\cdot)$ for $1 \le i \le n+1$.*

## 4. Proof of Theorem 1.2 : reduction via limiting arguments

We need to prove that: if $s > \frac{nd}{2}$ and

$$\Big(\frac{1}{p_1}, \cdots, \frac{1}{p_n}\Big) \in \mathrm{B}_n\Big(\frac{s}{d}\Big) := \Big\{(x_1, \cdots, x_n) \in (0,\infty)^n : \sum_{i=1}^n \max(x_i, 1/2) < \frac{s}{d}\Big\},$$

then we have

$$\big\|\mathrm{T}_{\mathbf{m}}(f_1, \ldots, f_n)\|_{L^p(\mathbb{R}^d)} \le C_{s,\delta}\, \|\mathbf{m}\|_{\mathscr{L}^2_{s,\delta}} \prod_{i=1}^n \|f_i\|_{h^{p_i}(\mathbb{R}^d)}.$$

We begin by replacing the multiplier $\mathbf{m}(x, \vec{\xi})$ with

$$\mathbf{m}_\lambda(x, \vec{\xi}) := \gamma(\lambda x)\mathbf{m}(x, \vec{\xi}), \quad 0 < \lambda \le 1;$$

here $\gamma$ is a fixed non-negative smooth function of compact support, with $\gamma(0) = 1$. For $f_i \in \mathscr{S}(\mathbb{R}^d)$, $1 \le i \le n$, we have

$$\mathrm{T}_{\mathbf{m}_\lambda}(f_1, \cdots, f_n)(x) = \gamma(\lambda x)\mathrm{T}_{\mathbf{m}}(f_1, \cdots, f_n)(x).$$

Therefore if we have

$$\big\|\mathrm{T}_{\mathbf{m}_\lambda}(f_1, \ldots, f_n)\big\|_{L^p(\mathbb{R}^d)} \le C\, \|\mathbf{m}\|_{\mathscr{L}^2_{s,\delta}} \prod_{i=1}^n \|f_i\|_{h^{p_i}(\mathbb{R}^d)},$$

uniformly in $0 < \lambda \le 1$. Then since $|\mathrm{T}_{\mathbf{m}_\lambda}(f_1, \ldots, f_n)(x)| \nearrow |\mathrm{T}_{\mathbf{m}}(f_1, \ldots, f_n)(x)|$ as $\lambda \to 0$, by the monotone convergence theorem we get

$$\big\|\mathrm{T}_{\mathbf{m}}(f_1, \ldots, f_n)\big\|_{L^p(\mathbb{R}^d)} \le C\, \|\mathbf{m}\|_{\mathscr{L}^2_{s,\delta}} \prod_{i=1}^n \|f_i\|_{h^{p_i}(\mathbb{R}^d)}.$$

By applying Lemma 3.11, we express $\mathrm{T}_{\mathbf{m}_\lambda}(f_1, f_2, \cdots, f_n)$ as a finite sum of the form

$$\mathscr{T}_{\mathbf{m}_\lambda}(f_1, \cdots, f_n)(x) := \sum_{j=-1}^{\infty} \mathscr{T}^j_{\mathbf{m}_\lambda}(f_1, \cdots, f_n)(x)$$

where
(4.1)

$$\mathscr{T}^{-1}_{\mathbf{m}_\lambda}(f_1, \cdots, f_n)(x) := \gamma(\lambda x)\int_{(\mathbb{R}^d)^n} \mathbf{m}(x, \vec{\xi})\widehat{\Phi}(\vec{\xi})\Big(\prod_{i=1}^n \widehat{f_i}(\xi_i)\Big)e^{2\pi i\sum_{i=1}^n\langle x,\xi_i\rangle}d\vec{\xi},$$

$$\mathscr{T}^j_{\mathbf{m}_\lambda}(f_1, \cdots, f_n)(x) := \gamma(\lambda x)\int_{(\mathbb{R}^d)^n} \mathbf{m}(x, \vec{\xi})\widehat{\Psi}(2^{-j}\vec{\xi})\Big(\prod_{i=1}^n \big(\widehat{\Phi}^i_j(\xi_i)\widehat{f_i}(\xi_i)\big)\Big)\widehat{\Phi}^{n+1}_j(\xi_{n+1})e^{2\pi i\sum_{i=1}^n\langle x,\xi_i\rangle}d\vec{\xi}$$

for $j \ge 0$, where $\vec{\xi} = (\xi_1, \cdots, \xi_n)$ and $\xi_{n+1} = -(\xi_1 + \cdots + \xi_n)$.

First we treat the term $\sum_{j\ge 0} \mathscr{T}^j_{\mathbf{m}_\lambda}(f_1, \cdots, f_n)(x)$. Since $\mathbf{m}(x, \vec{\xi})$ is a bounded function, we have

$$|\mathscr{T}^j_{\mathbf{m}_\lambda}(f_1, \cdots, f_n)(x)| \lesssim \gamma(\lambda x) \int_{(\mathbb{R}^d)^n} \big|\widehat{\Psi}(2^{-j}\vec{\xi})\big| \prod_{i=1}^{n} |\widehat{f}_i(\xi_i)|\, d\vec{\xi}.$$

Since $\sum_{j\ge 0} |\widehat{\Psi}(2^{-j}\vec{\xi})| \lesssim 1$ and $f_1, \ldots, f_n \in \mathscr{S}(\mathbb{R}^d)$, this implies that

$$\sum_{j=0}^{\infty} \|\mathscr{T}^j_{\mathbf{m}_\lambda}(f_1, \cdots, f_n)\|_{L^1(\mathbb{R}^d)} < \infty, \quad \text{and} \quad \sum_{j=0}^{\infty} \|\mathscr{T}^j_{\mathbf{m}_\lambda}(f_1, \cdots, f_n)\|_{L^2(\mathbb{R}^d)} < \infty. \tag{4.2}$$

By (4.2) we have $\sum_{j=0}^{N_1} \mathscr{T}^j_{\mathbf{m}_\lambda}(f_1, \cdots, f_n) \to \sum_{j=0}^{\infty} \mathscr{T}^j_{\mathbf{m}_\lambda}(f_1, \cdots, f_n)$ in $L^2(\mathbb{R}^d)$ as $N_1 \to \infty$. Thus by Lemma 3.10 it suffices to prove that

$$\Big\| \sum_{j=0}^{N_1} \mathscr{T}^j_{\mathbf{m}_\lambda}(f_1, \cdots, f_n) \Big\|_{L^p(\mathbb{R}^d)} \le C\, \|\mathbf{m}\|_{\mathscr{L}^2_{s,\delta}} \prod_{i=1}^{n} \|f_i\|_{h^{p_i}(\mathbb{R}^d)},$$

uniformly in $0 < \lambda \le 1$ and $N_1$. Fix $N_1$. Let $\varphi$ and $\psi$ be as in (3.13), then

$$\begin{aligned}
&\mathscr{F}\Big(\sum_{j=0}^{N_1} \mathscr{T}^j_{\mathbf{m}_\lambda}(f_1, \cdots, f_n)\Big)(\eta) \\
&= \sum_{j=0}^{N_1} \mathscr{F}\big(\mathscr{T}^j_{\mathbf{m}_\lambda}(f_1, \cdots, f_n)\big)(\eta)\Big(\widehat{\varphi}_{j-10}(\eta) + \sum_{k=j-9}^{j+9} \widehat{\psi}_k(\eta) + \sum_{k=j+10}^{\infty} \widehat{\psi}_k(\eta)\Big)
\end{aligned}$$

where $\mathscr{F}(f)$ denotes the Fourier transform of $f$. Since

$$\begin{aligned}
&\left(\sum_{j=0}^{N_1} \mathscr{F}\big(\mathscr{T}^j_{\mathbf{m}_\lambda}(f_1, \cdots, f_n)\big)(\eta)\Big(\widehat{\varphi}_{j-10}(\eta) + \sum_{k=j-9}^{j+9} \widehat{\psi}_k(\eta) + \sum_{k=j+10}^{N_2} \widehat{\psi}_k(\eta)\Big)\right) \\
&\to \mathscr{F}\Big(\sum_{j=0}^{N_1} \mathscr{T}^j_{\mathbf{m}_\lambda}(f_1, \cdots, f_n)\Big)(\eta)
\end{aligned}$$

in $L^2(\mathbb{R}^d)$ as $N_2 \to \infty$, by Plancherel's Theorem

$$\begin{aligned}
&\mathscr{F}^{-1}\left(\sum_{j=0}^{N_1} \mathscr{F}\big(\mathscr{T}^j_{\mathbf{m}_\lambda}(f_1, \cdots, f_n)\big)(\eta)\Big(\widehat{\varphi}_{j-10}(\eta) + \sum_{k=j-9}^{j+9} \widehat{\psi}_k(\eta) + \sum_{k=j+10}^{N_2} \widehat{\psi}_k(\eta)\Big)\right)(x) \\
&\to \sum_{j=0}^{N_1} \mathscr{T}^j_{\mathbf{m}_\lambda}(f_1, \cdots, f_n)(x)
\end{aligned}$$

in $L^2(\mathbb{R}^d)$ as $N_2 \to \infty$. Therefore by Lemma 3.10 it suffices to prove that

$$\Big\| \mathrm{I}^{N_1,N_2}_{\mathbf{m}_\lambda}(f_1, \cdots, f_n) \Big\|_{L^p(\mathbb{R}^d)} \le C\, \|\mathbf{m}\|_{\mathscr{L}^2_{s,\delta}} \prod_{i=1}^{n} \|f_i\|_{h^{p_i}(\mathbb{R}^d)},$$

$$\Big\| \mathrm{II}^{N_1}_{\mathbf{m}_\lambda}(f_1, \cdots, f_n) \Big\|_{L^p(\mathbb{R}^d)} \le C\, \|\mathbf{m}\|_{\mathscr{L}^2_{s,\delta}} \prod_{i=1}^{n} \|f_i\|_{h^{p_i}(\mathbb{R}^d)},$$

$$\Big\| \mathrm{III}^{N_1}_{\mathbf{m}_\lambda}(f_1, \cdots, f_n) \Big\|_{L^p(\mathbb{R}^d)} \le C\, \|\mathbf{m}\|_{\mathscr{L}^2_{s,\delta}} \prod_{i=1}^{n} \|f_i\|_{h^{p_i}(\mathbb{R}^d)},$$

uniformly in $0 < \lambda \le 1$, $N_1$ and $N_2$ with $N_1 + 100 \le N_2$, where

$$
\begin{aligned}
\mathrm{I}^{N_1,N_2}_{\mathbf{m}_\lambda}(f_1,\cdots,f_n)(x) &:= \sum_{j=0}^{N_1}\sum_{k=j+10}^{N_2} \mathscr{T}^j_{\mathbf{m}_\lambda}(f_1,\cdots,f_n)*\psi_k(x),\\
\mathrm{II}^{N_1}_{\mathbf{m}_\lambda}(f_1,\cdots,f_n)(x) &:= \sum_{j=0}^{N_1} \mathscr{T}^j_{\mathbf{m}_\lambda}(f_1,\cdots,f_n)*\varphi_{j-10}(x),\\
\mathrm{III}^{N_1}_{\mathbf{m}_\lambda}(f_1,\cdots,f_n)(x) &:= \sum_{j=0}^{N_1}\sum_{k=j-9}^{j+9} \mathscr{T}^j_{\mathbf{m}_\lambda}(f_1,\cdots,f_n)*\psi_k(x).
\end{aligned} \tag{4.3}
$$

The estimates for the term $\mathscr{T}^{-1}_{\mathbf{m}_\lambda}(f_1,\cdots,f_n)$ in (4.1) are similar to those for $\mathrm{I}^{0,N_2}_{\mathbf{m}_\lambda}(f_1,\cdots,f_n)$, and $\mathrm{II}^{0}_{\mathbf{m}_\lambda}(f_1,\cdots,f_n)$ in (4.3). The proof for $\mathscr{T}^{-1}_{\mathbf{m}_\lambda}(f_1,\cdots,f_n)$ will be sketched briefly in Section 8.

**Remark 5.** For $\mathrm{I}^{N_1,N_2}_{\mathbf{m}_\lambda}(f_1,\cdots,f_n)$ and $\mathrm{II}^{N_1}_{\mathbf{m}_\lambda}(f_1,\cdots,f_n)$, we make use of the derivative condition in (1.20) concerning the space variable $x$ to obtain the summability over indices $j$ and $k$ above.

## 5. Proof of Theorem 1.2 : Estimates for the term $\mathrm{III}^{N_1}_{\mathbf{m}_\lambda}(f_1,\cdots,f_n)$ in (4.3)

By (4.3), we have

$$
\mathrm{III}^{N_1}_{\mathbf{m}_\lambda}(f_1,\cdots,f_n)(x) = \sum_{i=-9}^{9}\sum_{j=0}^{N_1} \mathscr{T}^j_{\mathbf{m}_\lambda}(f_1,\cdots,f_n)*\psi_{j+i}(x).
$$

Then by Lemma 3.2, we have

$$
\begin{aligned}
\|\mathrm{III}^{N_1}_{\mathbf{m}_\lambda}(f_1,\cdots,f_n)\|_{L^p} &\le \|\mathrm{III}^{N_1}_{\mathbf{m}_\lambda}(f_1,\cdots,f_n)\|_{h^p}\\
&\lesssim \sum_{i=-9}^{9}\Big\|\Big(\sum_{j=0}^{N_1}\Big|\mathscr{T}^j_{\mathbf{m}_\lambda}(f_1,\cdots,f_n)*\psi_{j+i}\Big|^2\Big)^{1/2}\Big\|_{L^p}.
\end{aligned}
$$

We will only consider the case $i = 0$ in the previous summation. Note that

$$
\begin{aligned}
&\mathscr{T}^j_{\mathbf{m}_\lambda}(f_1,\cdots,f_n)*\psi_j(x)\\
&= \int_{\mathbb{R}^{(n+1)d}} \mathbf{m}_\lambda(x-y,\vec{\xi})\widehat{\Psi}(2^{-j}\vec{\xi})\,\widehat{\Phi}^{n+1}_j(-\xi_1-\cdots-\xi_n)\\
&\quad\times\Big(\prod_{i=1}^n\big(\widehat{\Phi}^i_j(\xi_i)\widehat{f}_i(\xi_i)\big)\Big)e^{2\pi i\langle x-y,\xi_1+\cdots+\xi_n\rangle}\Big(\int_{\mathbb{R}^d}\widehat{\psi_j}(\eta)e^{2\pi i\langle y,\eta\rangle}d\eta\Big)dy\,d\vec{\xi}.
\end{aligned}
$$

Let

$$
\mathbf{m}_j(x,\vec{\xi}) := \mathbf{m}(x,2^j\vec{\xi})\widehat{\Psi}(\vec{\xi})\widehat{\Phi}^{n+1}(-\xi_1-\cdots-\xi_n).
$$

Then by using the identity

$$
\mathbf{m}(x,\vec{\xi})\widehat{\Psi}(2^{-j}\vec{\xi})\widehat{\Phi}^{n+1}_j(-\xi_1-\cdots-\xi_n) = \int_{(\mathbb{R}^d)^n}\mathscr{F}^{-1}[\mathbf{m}_j(x,\,\vec{\cdot}\,)](\vec{z})e^{-2\pi i\langle\vec{z},2^{-j}\vec{\xi}\rangle}d\vec{z},
$$

we obtain that

$$
\begin{aligned}
&\mathscr{T}^j_{\mathbf{m}_\lambda}(f_1,\cdots,f_n)*\psi_j(x)\\
&= \int_{(\mathbb{R}^d)^{n+1}}\gamma(\lambda(x-y))\mathscr{F}^{-1}[\mathbf{m}_j(x-y,\,\vec{\cdot}\,)](\vec{z})\Big(\prod_{i=1}^n f_i*\Phi^i_j(x-y-2^{-j}z_i)\Big)\psi_j(y)\,d\vec{z}\,dy.
\end{aligned} \tag{5.1}
$$

**Lemma 5.1.** *Let* $\mathbf{m}_j(x,\vec{\xi}) := \mathbf{m}(x,2^j\vec{\xi})\widehat{\Psi}(\vec{\xi})\widehat{\Phi}^{n+1}(-\xi_1-\cdots-\xi_n)$*. Then for* $s \geq 0$

$$\int |\mathscr{F}^{-1}[\mathbf{m}_j(x,\vec{\cdot})](\vec{z})|^2(1+|\vec{z}|)^{2s}d\vec{z} \lesssim \int |\mathscr{F}^{-1}[\mathbf{m}(x,2^j\vec{\cdot})\widehat{\Psi}(\cdot)](\vec{z})|^2(1+|\vec{z}|)^{2s}d\vec{z}.$$

*Proof of Lemma 5.1.* By adopting smooth bump function $\Psi'$ we write

$$\mathbf{m}_j(x,\vec{\xi}) = \mathbf{m}(x,2^j\vec{\xi})\widehat{\Psi}(\vec{\xi})\widehat{\Psi'}(\vec{\xi})\widehat{\Phi}^{n+1}(-\xi_1-\cdots-\xi_n).$$

Then by using

$$\left|\mathscr{F}^{-1}\left[\widehat{\Psi'}(\vec{\xi})\widehat{\Phi}^{n+1}(-\xi_1-\cdots-\xi_n)\right](\vec{z})\right| \leq C_N(1+|\vec{z}|)^{-N} \quad \text{for any } N>0,$$

and Hölder's inequality we have

$$\begin{aligned}
\left|\mathscr{F}^{-1}[\mathbf{m}_j(x,\vec{\cdot})](\vec{z})\right|^2 &\lesssim \left|\int_{\mathbb{R}^{dn}} |\mathscr{F}^{-1}[\mathbf{m}(x,2^j\vec{\cdot})\widehat{\Psi}(\vec{\cdot})](\vec{z}-\vec{y})|(1+|\vec{y}|)^{-N}d\vec{y}\right|^2 \\
&\lesssim \int_{\mathbb{R}^{dn}} \left|\mathscr{F}^{-1}[\mathbf{m}(x,2^j\vec{\cdot})\widehat{\Psi}(\vec{\cdot})](\vec{z}-\vec{y})\right|^2(1+|\vec{y}|)^{-N}d\vec{y}.
\end{aligned} \tag{5.2}$$

If $N > nd+2s$, then by considering the following integral into three cases $|\vec{z}| > 2|\vec{y}|$, $|\vec{z}| < 2|\vec{y}|$, or $|\vec{z}| \sim |\vec{y}|$, we have

$$\sup_{\vec{z}} \int \frac{(1+|\vec{z}+\vec{y}|)^{2s}}{(1+|\vec{z}|)^{2s}} \frac{1}{(1+|\vec{y}|)^N} d\vec{y} \lesssim 1. \tag{5.3}$$

Now by (5.2) and (5.3)

$$\begin{aligned}
&\int |\mathscr{F}^{-1}[\mathbf{m}_j(x,\vec{\cdot})](\vec{z})|^2(1+|\vec{z}|)^{2s}d\vec{z} \\
&\lesssim \iint \left|\mathscr{F}^{-1}[\mathbf{m}(x,2^j\vec{\cdot})\widehat{\Psi}(\vec{\cdot})](\vec{z}-\vec{y})\right|^2 \frac{(1+|\vec{z}|)^{2s}}{(1+|\vec{y}|)^N} d\vec{z}d\vec{y} \\
&= \iint \left|\mathscr{F}^{-1}[\mathbf{m}(x,2^j\vec{\cdot})\widehat{\Psi}(\vec{\cdot})](\vec{z})\right|^2 \frac{(1+|\vec{z}+\vec{y}|)^{2s}}{(1+|\vec{y}|)^N} d\vec{z}d\vec{y} \\
&= \int \left|\mathscr{F}^{-1}[\mathbf{m}(x,2^j\vec{\cdot})\widehat{\Psi}(\vec{\cdot})](\vec{z})\right|^2(1+|\vec{z}|)^{2s}\left(\int \frac{(1+|\vec{z}+\vec{y}|)^{2s}}{(1+|\vec{z}|)^{2s}} \frac{1}{(1+|\vec{y}|)^N} d\vec{y}\right)d\vec{z} \\
&\lesssim \int \left|\mathscr{F}^{-1}[\mathbf{m}(x,2^j\vec{\cdot})\widehat{\Psi}(\vec{\cdot})](\vec{z})\right|^2(1+|\vec{z}|)^{2s}d\vec{z}.
\end{aligned}$$

□

Let

$$\Gamma_0 := \left\{\vec{z} \in (\mathbb{R}^d)^n : |\vec{z}| \leq 1\right\}, \quad \Gamma_M := \left\{\vec{z} \in (\mathbb{R}^d)^n : 2^{M-1} < |\vec{z}| \leq 2^M\right\}, \quad M \geq 1.$$

Then by (5.1) we have

$$|\mathscr{T}^j_{\mathbf{m}_\lambda}(f_1,\cdots,f_n)*\psi_j(x)| \leq \sum_{M\geq 0} \mathrm{C}^{M,j}_{\mathbf{m}_\lambda}(f_1,\cdots,f_n)(x) \tag{5.4}$$

where

$$\begin{aligned}
\mathrm{C}^{M,j}_{\mathbf{m}_\lambda}(f_1,\cdots,f_n)(x) := &\int_{\mathbb{R}^d}\int_{\Gamma_M} \left|\gamma(\lambda(x-y))\mathscr{F}^{-1}[\mathbf{m}_j(x-y,\vec{\cdot})](\vec{z})\right| \\
&\times\Big(\prod_{i=1}^n |f_i*\Phi^i_j(x-y-2^{-j}z_i)|\Big)|\psi_j(y)|\,d\vec{z}dy.
\end{aligned} \tag{5.5}$$

Since the case $M=0$ is similar to the case $M\geq 1$, we only consider the case $M\geq 1$. By (5.5), we have

$$\mathrm{C}^{M,j}_{\mathbf{m}_\lambda}(f_1,\cdots,f_n)(x)\lesssim 2^{Mnd}\int_{\mathbb{R}^d}\int_{|\vec{\mathbf{z}}|\sim 1}\left|\gamma(\lambda(x-y))\mathscr{F}^{-1}[\mathbf{m}_j(x-y,\vec{\cdot}\,)](2^M\vec{\mathbf{z}})\right| \times\Big(\prod_{i=1}^n|f_i*\Phi^i_j(x-y-2^{-j+M}z_i)|\Big)|\psi_j(y)|\,d\vec{\mathbf{z}}dy.$$

Then by using the Hölder's inequality with $\vec{\mathbf{z}}$ variable and Lemma 5.1

$$\mathrm{C}^{M,j}_{\mathbf{m}_\lambda}(f_1,\cdots,f_n)(x)\lesssim 2^{-Ms+\frac{Mnd}{2}}\Big(\sup_{x\in\mathbb{R}^d}\big\|\mathbf{m}(x,2^j\vec{\cdot}\,)\widehat{\Psi}(\vec{\cdot}\,)\big\|_{L^2_s((\mathbb{R}^d)^n)}\Big) \times\int_{\mathbb{R}^d}\Big(\prod_{i=1}^n\Big(\int_{|z_i|\lesssim 1}\big|f_i*\Phi^i_j(x-y-2^{-j+M}z_i)\big|^2dz_i\Big)^{1/2}\Big)|\psi_j(y)|dy. \tag{5.6}$$

Let $|y|\sim 2^{-j+M+l}$ for some $l\geq 0$, then $y=2^{-j+M+l}y'$ for some $|y'|\lesssim 1$. By the change of variables $y'+2^{-l}z_i\to z_i$ we have

$$\begin{aligned}\Big(\int_{|z_i|\lesssim 1}\big|f_i*\Phi^i_j(x-y-2^{-j+M}z_i)\big|^2dz_i\Big)^{1/2}&=\Big(\int_{|z_i|\lesssim 1}\big|f_i*\Phi^i_j(x-2^{-j+M+l}(y'+2^{-l}z_i))\big|^2dz_i\Big)^{1/2}\\&=2^{\frac{ld}{2}}\Big(\int_{|z_i|\lesssim 1}\big|f_i*\Phi^i_j(x-2^{-j+M+l}z_i)\big|^2dz_i\Big)^{1/2}.\end{aligned}$$

Then by Lemma 3.7, for any $0<q_i\leq 2$,

$$\Big(\int_{|z_i|\lesssim 1}\big|f_i*\Phi^i_j(x-y-2^{-j+M}z_i)\big|^2dz_i\Big)^{1/2}\lesssim 2^{\frac{ld}{2}}2^{(M+l)d(\frac{1}{q_i}-\frac12)}\mathrm{M}_{q_i}(f_i*\Phi^i_j)(x). \tag{5.7}$$

Therefore by (5.7), if we take $N$ large enough, then

$$\begin{aligned}&\int\Big(\prod_{i=1}^n\Big(\int_{|z_i|\lesssim 1}\big|f_i*\Phi^i_j(x-y-2^{-j+M}z_i)\big|^2dz_i\Big)^{1/2}\Big)|\psi_j(y)|dy\\&\lesssim\int_{|y|\lesssim 2^{-j+M}}\Big(\cdots\Big)dy+\sum_{l\geq 0}\int_{|y|\sim 2^{-j+M+l}}\Big(\cdots\Big)dy\\&\lesssim\int_{|y|\lesssim 2^{-j+M}}\Big(\prod_{i=1}^n 2^{Md(\frac{1}{q_i}-\frac12)}\mathrm{M}_{q_i}(f_i*\Phi^i_j)(x)\Big)|\psi_j(y)|dy\\&\quad+\sum_{l\geq 0}\prod_{i=1}^n\Big(2^{\frac{ld}{2}}2^{(M+l)d(\frac{1}{q_i}-\frac12)}\mathrm{M}_{q_i}(f_i*\Phi^i_j)(x)\Big)\Big(\frac{1}{(1+2^{M+l})^N}\Big)\\&\lesssim\prod_{i=1}^n\Big(2^{Md(\frac{1}{q_i}-\frac12)}\mathrm{M}_{q_i}(f_i*\Phi^i_j)(x)\Big),\end{aligned} \tag{5.8}$$

where we use $|\psi_j(y)|\leq C_N2^{jd}(1+|2^jy|)^{-N-d-1}$ for the second inequality.

By (5.6), (5.8), and Lemma 5.1, for $0<q_j\leq 2$ we have

$$\begin{aligned}&\mathrm{C}^{M,j}_{\mathbf{m}_\lambda}(f_1,\cdots,f_n)(x)\\&\lesssim\Big(2^{-Ms+\frac{Mnd}{2}+\sum_{i=1}^n Md\big(\frac{1}{q_i}-\frac12\big)}\Big)\Big(\sup_{x\in\mathbb{R}^d}\big\|\mathbf{m}(x,2^j\vec{\cdot}\,)\widehat{\Psi}(\vec{\cdot}\,)\big\|_{L^2_s((\mathbb{R}^d)^n)}\Big)\Big(\prod_{i=1}^n\mathrm{M}_{q_i}(f_i*\Phi^i_j)(x)\Big).\end{aligned}$$

By Lemma 3.11, there exists at one $\widehat{\Phi}^i$ $(1 \le i \le n)$ that is compactly supported away from the origin. Without loss of generality, let $\widehat{\Phi}^1$ be a such function. Then

$$\begin{aligned}\mathrm{C}^{M,j}_{\mathbf{m}_\lambda}(f_1,\cdots,f_n)(x) &\lesssim \Big(2^{-Ms+\frac{Mnd}{2}+\sum_{i=1}^n Md\left(\frac{1}{q_i}-\frac{1}{2}\right)}\Big)\Big(\sup_{j\ge 0}\sup_{x\in\mathbb{R}^d}\big\|\mathbf{m}(x,2^j\vec{\cdot})\widehat{\Psi}(\vec{\cdot})\big\|_{L^2_s((\mathbb{R}^d)^n)}\Big)\\ &\quad\times \mathrm{M}_{q_1}(f_1*\Phi^1_j)(x)\Big(\prod_{i=2}^n \mathrm{M}_{q_i}(\sup_{j\ge 0}|f_i*\Phi^i_j|)(x)\Big).\end{aligned} \tag{5.9}$$

By (5.4) and (5.9) we have

$$\begin{aligned}&\Big(\sum_{j=0}^{N_1}|\mathscr{T}^j_{\mathbf{m}_\lambda}(f_1,\cdots,f_n)*\psi_j(x)|^2\Big)^{1/2}\\ &\le \sum_{M\ge 0}\Big(\sum_{j=0}^{N_1}\big|\mathrm{C}^{M,j}_{\mathbf{m}_\lambda}(f_1,\cdots,f_n)(x)\big|^2\Big)^{1/2}\\ &\lesssim \sum_{M\ge 0}\Big(2^{-Ms+\frac{Mnd}{2}+\sum_{i=1}^n Md\left(\frac{1}{q_i}-\frac{1}{2}\right)}\Big)\Big(\sup_{j\ge 0}\sup_{x\in\mathbb{R}^d}\big\|\mathbf{m}(x,2^j\vec{\cdot})\widehat{\Psi}(\vec{\cdot})\big\|_{L^2_s((\mathbb{R}^d)^n)}\Big)\\ &\times\Big(\sum_{j=0}^{N_1}|\mathrm{M}_{q_1}(f_1*\Phi^1_j)(x)|^2\Big)^{\frac{1}{2}}\Big(\prod_{i=2}^n \mathrm{M}_{q_i}(\sup_{j\ge 0}|f_i*\Phi^i_j|)(x)\Big).\end{aligned} \tag{5.10}$$

Since $\frac{p}{p_1}+\cdots+\frac{p}{p_n}=1$, by Hölder's inequality, if $0<q_i<\min(2,p_i)$, then by Lemma 3.4 we have

$$\begin{aligned}&\Big\|\Big(\sum_{j=0}^{N_1}|\mathrm{M}_{q_1}(f_1*\Phi^1_j)|^2\Big)^{\frac{1}{2}}\Big(\prod_{i=2}^n \mathrm{M}_{q_i}(\sup_{j\ge 0}|f_i*\Phi^i_j|)\Big)\Big\|_{L^p}\\ &\lesssim \Big\|\Big(\sum_{j=0}^{N_1}|\mathrm{M}_{q_1}(f_1*\Phi^1_j)|^2\Big)^{\frac{1}{2}}\Big\|_{L^{p_1}}\prod_{i=2}^n\|\mathrm{M}_{q_i}(\sup_{j\ge 0}|f_i*\Phi^i_j|)\|_{L^{p_i}}\\ &\lesssim \Big\|\Big(\sum_{j=0}^{N_1}|f_1*\Phi^1_j|^2\Big)^{\frac{1}{2}}\Big\|_{L^{p_1}}\prod_{i=2}^n\|\sup_{j\ge 0}|f_i*\Phi^i_j|\|_{L^{p_i}}\\ &\lesssim \prod_{i=1}^n\|f_i\|_{h^{p_i}}.\end{aligned} \tag{5.11}$$

Then by (5.10) and (5.11), if $0<q_i<\min(2,p_i)$ for all $1\le i\le n$, then

$$\begin{aligned}\|\mathrm{III}^{N_1}_{\mathbf{m}_\lambda}(f_1,\cdots,f_n)\|_{L^p} &\lesssim \sum_{i=-9}^{9}\Big\|\Big(\sum_{j=0}^{N_1}|\mathscr{T}^j_{\mathbf{m}_\lambda}(f_1,\cdots,f_n)*\psi_{j+i}|^2\Big)^{1/2}\Big\|_{L^p}\\ &\lesssim \Big(\sup_{j\ge 0}\sup_{x\in\mathbb{R}^d}\big\|\mathbf{m}(x,2^j\vec{\cdot})\widehat{\Psi}(\vec{\cdot})\big\|_{L^2_s((\mathbb{R}^d)^n)}\Big)\prod_{i=1}^n\|f_i\|_{h^{p_i}},\end{aligned} \tag{5.12}$$

uniformly in $0<\lambda\le 1$ and $N_1$ when $\frac{s}{d}>\sum_{i=1}^n\frac{1}{q_i}$. Therefore, by taking $q_i\nearrow\min(2,p_i)$ for $1\le i\le n$, we have (5.12) when

$$\Big(\frac{1}{p_1},\cdots,\frac{1}{p_n}\Big)\in \mathrm{B}_n\Big(\frac{s}{d}\Big)=\Big\{(x_1,\cdots,x_n)\in(0,\infty)^n : \sum_{i=1}^n\max(x_i,1/2)<\frac{s}{d}\Big\}.$$

## 6. Proof of Theorem 1.2 : Estimates for the term $\mathrm{I}^{N_1,N_2}_{\mathbf{m}_\lambda}(f_1,\cdots,f_n)$ in (4.3)

By (4.3) we have

$$\mathrm{I}^{N_1,N_2}_{\mathbf{m}_\lambda}(f_1,\cdots,f_n)(x)=\sum_{j=0}^{N_1}\sum_{k=j+10}^{N_2}\mathscr{T}^j_{\mathbf{m}_\lambda}(f_1,\cdots,f_n)*\psi_k(x).$$

Note that

$$\begin{aligned}&\mathscr{T}^j_{\mathbf{m}_\lambda}(f_1,\cdots,f_n)*\psi_k(x)\\&=\int_{\mathbb{R}^{(n+1)d}}\mathbf{m}_\lambda(x-y,\vec{\xi})\widehat{\Psi}(2^{-j}\vec{\xi})\,\widehat{\Phi}^{n+1}_j(-\xi_1-\cdots-\xi_n)\\&\quad\times\Big(\prod_{i=1}^n\big(\widehat{\Phi}^i_j(\xi_i)\widehat{f_i}(\xi_i)\big)\Big)e^{2\pi i\langle x-y,\xi_1+\cdots+\xi_n\rangle}\Big(\int_{\mathbb{R}^d}\widehat{\psi_k}(\eta)e^{2\pi i\langle y,\eta\rangle}d\eta\Big)dyd\vec{\xi}.\end{aligned}$$

If $k\ge j+10$, then $\widehat{\Phi}^{n+1}_j(-\xi_1-\cdots-\xi_n)\widehat{\psi_k}(\eta)\neq 0$ only if $|\xi_1+\cdots+\xi_n-\eta|\sim 2^k$. The term $\widehat{\Phi}^{n+1}_j(-\xi_1-\cdots-\xi_n)\widehat{\psi_k}(\eta)$ can be written as a finite sum of the form:

$$\widehat{\Phi}^{n+1}_j(-\xi_1-\cdots-\xi_n)\widehat{\psi_k}(\eta)\prod_{i=1}^d\phi^i\Big(\frac{\xi_1^i+\cdots+\xi_n^i-\eta^i}{2^k}\Big)$$

where $\phi^i(t)$ are smooth functions that are supported in $|t|\lesssim 1$ and at least one of $\phi^i(t)$ is supported in $|t|\sim 1$. Thus $\mathscr{T}^j_{\mathbf{m}_\lambda}(f_1,\cdots,f_n)*\psi_k(x)$ can be written as a finite sum of the form:

$$\begin{aligned}&\mathscr{T}^j_{\mathbf{m}_\lambda}(f_1,\cdots,f_n)*\psi_k(x)\\&:=\int_{\mathbb{R}^{(n+2)d}}\prod_{i=1}^d\phi^i\Big(\frac{\xi_1^i+\cdots+\xi_n^i-\eta^i}{2^k}\Big)\mathbf{m}_\lambda(x-y,\vec{\xi})\widehat{\Psi}(2^{-j}\vec{\xi})\,\widehat{\Phi}^{n+1}_j(-\xi_1-\cdots-\xi_n)\\&\quad\times\Big(\prod_{i=1}^n\big(\widehat{\Phi}^i_j(\xi_i)\widehat{f_i}(\xi_i)\big)\Big)\widehat{\psi_k}(\eta)e^{2\pi i\langle y,\eta-(\xi_1+\cdots+\xi_n)\rangle}e^{2\pi i\langle x,\xi_1+\cdots+\xi_n\rangle}dyd\eta d\vec{\xi}.\end{aligned}$$

where $\phi^l(t)$ is supported in $|t|\sim 1$ for some $1\le l\le d$. Then by using integration by parts via

$$\partial_{y_l}e^{2\pi i\langle y,\eta-(\xi_1+\cdots+\xi_n)\rangle}=2\pi i(\eta^l-\xi_1^l-\cdots-\xi_n^l)e^{2\pi i\langle y,\eta-(\xi_1+\cdots+\xi_n)\rangle},$$

we obtain that

$$\begin{aligned}&\mathscr{T}^j_{\mathbf{m}_\lambda}(f_1,\cdots,f_n)*\psi_k(x)\\&=-\int_{\mathbb{R}^{(n+2)d}}\frac{e^{2\pi i\langle x,\xi_1+\cdots+\xi_n\rangle}}{2\pi i(\eta^l-\xi_1^l-\cdots-\xi_n^l)}\prod_{i=1}^d\phi^i\Big(\frac{\xi_1^i+\cdots+\xi_n^i-\eta^i}{2^k}\Big)\big[\partial_{y_l}\mathbf{m}_\lambda(x-y,\vec{\xi})\widehat{\Psi}(2^{-j}\vec{\xi})\big]\\&\quad\times\Big(\prod_{i=1}^n\big(\widehat{\Phi}^i_j(\xi_i)\widehat{f_i}(\xi_i)\big)\Big)\widehat{\Phi}^{n+1}_j(-\xi_1-\cdots-\xi_n)\widehat{\psi_k}(\eta)e^{2\pi i\langle y,\eta-(\xi_1+\cdots+\xi_n)\rangle}dyd\eta d\vec{\xi}.\end{aligned}$$

Then by using the identity

$$\partial_{y_l}\mathbf{m}_\lambda(x-y,\vec{\xi})\widehat{\Psi}(2^{-j}\vec{\xi})=\int_{(\mathbb{R}^d)^n}\mathscr{F}^{-1}\Big(\partial_{y_l}\big[\mathbf{m}_\lambda(x-y,2^j\vec{\cdot})\widehat{\Psi}(\vec{\cdot})\big]\Big)(\vec{z})e^{-2\pi i\langle\vec{z},2^{-j}\vec{\xi}\rangle}d\vec{z},$$

we have

$$
\begin{aligned}
&\mathcal{T}^j_{\mathbf{m}_\lambda}(f_1,\cdots,f_n)*\psi_k(x)\\
(6.1)\quad &= -\int_{\mathbb{R}^{(2n+2)d}} \frac{e^{2\pi i\langle(\vec{\xi},\eta),(x-y-2^{-j}z_1,\cdots,x-y-2^{-j}z_n,y)\rangle}}{2\pi i(\eta^l-\xi_1^l-\cdots-\xi_n^l)} \prod_{i=1}^d \phi^i\Big(\frac{\xi_1^i+\cdots+\xi_n^i-\eta^i}{2^k}\Big)\widehat{\Phi}_j^{n+1}(-\xi_1-\cdots-\xi_n)\\
&\quad \times \mathcal{F}^{-1}\big(\partial_{y_l}\big[\mathbf{m}_\lambda(x-y,2^j\vec{\cdot})\widehat{\Psi}(\vec{\cdot})\big]\big)(\vec{\mathbf{z}})\Big(\prod_{i=1}^n\big(\widehat{\Phi}_j^i(\xi_i)\widehat{f}_i(\xi_i)\big)\Big)\widehat{\psi}_k(\eta)dyd\eta d\vec{\mathbf{z}}d\vec{\boldsymbol{\xi}}.
\end{aligned}
$$

**Lemma 6.1.** *Let* $|j-k|\geq 10$, $|\xi_1^l+\cdots+\xi_n^l-\eta^l|\sim 2^{\max(k,j)}$, *and let*

$$
m_{k,j}(\vec{\xi},\eta) := \frac{\prod_{i=1}^d \phi^i\Big(\frac{\xi_1^i+\cdots+\xi_n^i-\eta^i}{2^{\max(k,j)}}\Big)}{2\pi i(\eta^l-\xi_1^l-\cdots-\xi_n^l)}\Big(\prod_{i=1}^n \widehat{\Phi}_j^i(\xi_i)\Big)\widehat{\Phi}_j^{n+1}(-\xi_1-\cdots-\xi_n)\widehat{\psi}_k(\eta).
$$

*Then for any positive integer* $N$ *we have*

$$
\begin{aligned}
&\Big|\int_{\mathbb{R}^{(n+1)d}} e^{2\pi i\langle(x_1,\cdots,x_n,x_{n+1}),(\xi_1,\cdots,\xi_n,\eta)\rangle} m_{k,j}(\xi_1,\cdots,\xi_n,\eta)\Big(\prod_{i=1}^n \widehat{f}_i(\xi_i)\Big)d\xi_1\cdots d\xi_n d\eta\Big|\\
&\leq C_N\frac{1}{\max(2^k,2^j)}\Big(\prod_{i=1}^n \omega_j^N*|f_i|(x_i)\Big)\omega_k^N(x_{n+1}),
\end{aligned}
$$

*where* $\omega_k^N(y)=\frac{2^{kd}}{(1+|2^ky|)^N}$.

*Proof of Lemma 6.1.* Note that

$$
\big|\partial^\beta_{\xi_i}(m_{k,j}(\xi_1,\cdots,\xi_n,\eta))\big| \lesssim \frac{1}{\max(2^k,2^j)}\frac{1}{2^{j|\beta|}},\quad \big|\partial^\beta_{\eta}(m_{k,j}(\xi_1,\cdots,\xi_n,\eta))\big| \lesssim \frac{1}{\max(2^k,2^j)}\frac{1}{2^{k|\beta|}},
$$

for all multi-indices $\beta$. The results follow from integration by parts via

$$
\begin{aligned}
\partial^\beta_{\xi_i}\big(e^{2\pi i\langle(x_1,\cdots,x_n,x_{n+1}),(\xi_1,\cdots,\xi_n,\eta)\rangle}\big) &= (2\pi i x_i)^\beta e^{2\pi i\langle(x_1,\cdots,x_n,x_{n+1}),(\xi_1,\cdots,\xi_n,\eta)\rangle},\\
\partial^\beta_{\eta}\big(e^{2\pi i\langle(x_1,\cdots,x_n,x_{n+1}),(\xi_1,\cdots,\xi_n,\eta)\rangle}\big) &= (2\pi i x_{n+1})^\beta e^{2\pi i\langle(x_1,\cdots,x_n,x_{n+1}),(\xi_1,\cdots,\xi_n,\eta)\rangle}.
\end{aligned}
$$

□

Since $k\geq j+10$, by (6.1) and Lemma 6.1, we have

$$
\begin{aligned}
|\mathcal{T}^j_{\mathbf{m}_\lambda}(f_1,\cdots,f_n)*\psi_k(x)| \leq C_N\frac{1}{2^k}\iint &\big|\mathcal{F}^{-1}\big(\partial_{y_l}\big[\mathbf{m}_\lambda(x-y,2^j\vec{\cdot})\widehat{\Psi}(\vec{\cdot})\big]\big)(\vec{\mathbf{z}})\big|\\
&\times\Big(\prod_{i=1}^n \omega_j^N*|f_i*\Phi_j^i|(x-y-2^{-j}z_i)\Big)(\omega_k^N(y))dyd\vec{\mathbf{z}}.
\end{aligned}
$$

Let

$$
\Gamma_0 := \big\{\vec{\mathbf{z}}\in(\mathbb{R}^d)^n : |\vec{\mathbf{z}}|\leq 1\big\},\quad \Gamma_M := \big\{\vec{\mathbf{z}}\in(\mathbb{R}^d)^n : 2^{M-1}<|\vec{\mathbf{z}}|\leq 2^M\big\},\quad M\geq 1.
$$

Then we have

$$
|\mathcal{T}^j_{\mathbf{m}_\lambda}(f_1,\cdots,f_n)*\psi_k(x)| \lesssim \sum_{M\geq 0} \mathrm{A}^{M,j,k}_{\mathbf{m}_\lambda}(f_1,\cdots,f_n)(x)
$$

where

$$\begin{aligned}\mathrm{A}^{M,j,k}_{\mathbf{m}_\lambda}(f_1,\cdots,f_n)(x) := \frac{1}{2^k}\int_{\mathbb{R}^d}\int_{\Gamma_M} \left|\mathscr{F}^{-1}\left(\partial_{y_l}\left[\mathbf{m}_\lambda(x-y,2^j\vec{\cdot})\widehat{\Psi}(\vec{\cdot})\right]\right)(\vec{\mathbf{z}})\right| \\ \times\left(\prod_{i=1}^n \omega_j^N * |f_i * \Phi_j^i|(x-y-2^{-j}z_i)\right)(\omega_k^N(y))d\vec{\mathbf{z}}dy.\end{aligned} \tag{6.2}$$

Since the case $M=0$ is similar to the case $M\geq 1$, we only consider the case $M\geq 1$. By (6.2) we have

$$\begin{aligned}\mathrm{A}^{M,j,k}_{\mathbf{m}_\lambda}(f_1,\cdots,f_n)(x) \lesssim \frac{2^{Mnd}}{2^k}\int_{\mathbb{R}^d}\int_{|\vec{\mathbf{z}}|\sim 1} \left|\mathscr{F}^{-1}\left(\partial_{y_l}\left[\mathbf{m}_\lambda(x-y,2^j\vec{\cdot})\widehat{\Psi}(\vec{\cdot})\right]\right)(2^M\vec{\mathbf{z}})\right| \\ \times\left(\prod_{i=1}^n \omega_j^N * |f_i * \Phi_j^i|(x-y-2^{-j+M}z_i)\right)(\omega_k^N(y))d\vec{\mathbf{z}}dy.\end{aligned}$$

By using the Hölder's inequality with the $\vec{\mathbf{z}}$ variable

$$\begin{aligned}|\mathrm{A}^{M,j,k}_{\mathbf{m}_\lambda}(f_1,\cdots,f_n)(x)| \lesssim \frac{2^{-Ms+\frac{Mnd}{2}}2^{j\delta}}{2^k}\sup_{x\in\mathbb{R}^d}\left(\sum_{|\alpha|\leq 1}2^{-j\delta}\left\|\partial_x^\alpha\mathbf{m}_\lambda(x,2^j\vec{\cdot})\widehat{\Psi}(\vec{\cdot})\right\|_{L^2_s((\mathbb{R}^d)^n)}\right) \\ \times\int\left(\prod_{i=1}^n\left(\int_{|z_i|\lesssim 1}\left|\omega_j^N * |f_i * \Phi_j^i|(x-y-2^{-j+M}z_i)\right|^2dz_i\right)^{1/2}\right)\omega_k^N(y)dy.\end{aligned} \tag{6.3}$$

Let $|y|\sim 2^{-j+M+l}$ for some $l\geq 0$, then $y=2^{-j+M+l}y'$ for some $|y'|\lesssim 1$. And by the change of variables $y'+2^{-l}z_i\to z_i$ we have

$$\begin{aligned}&\left(\int_{|z_i|\lesssim 1}\left|\omega_j^N * |f_i * \Phi_j^i|(x-y-2^{-j+M}z_i)\right|^2dz_i\right)^{1/2} \\ &=\left(\int_{|z_i|\lesssim 1}\left|\omega_j^N * |f_i * \Phi_j^i|(x-2^{-j+M+l}(y'+2^{-l}z_i))\right|^2dz_i\right)^{1/2} \\ &=2^{\frac{ld}{2}}\left(\int_{|z_i|\lesssim 1}\left|\omega_j^N * |f_i * \Phi_j^i|(x-2^{-j+M+l}z_i)\right|^2dz_i\right)^{1/2}.\end{aligned} \tag{6.4}$$

By Lemma 3.9, for any $0<q_i\leq 2$, if $N>(\frac{2}{q_i}+1)d$, then

$$\left(\int_{|z_i|\lesssim 1}\left|\omega_j^N * |f_i * \Phi_j^i|(x-2^{-j+M}z_i)\right|^2dz_i\right)^{1/2} \lesssim 2^{(M+l)d(\frac{1}{q_i}-\frac{1}{2})}\mathrm{M}_{q_i}(f_i * \Phi_j^i)(x). \tag{6.5}$$

Therefore by (6.4) and (6.5), if we take $N$ large enough, then

$$\begin{aligned}
&\int\Big(\prod_{i=1}^{n}\Big(\int_{|z_i|\lesssim 1}\big|\omega_j^N*|f_i*\Phi_j^i|(x-y-2^{-j+M}z_i)\big|^2dz_i\Big)^{1/2}\Big)\omega_k^N(y)dy\\
&\lesssim\int_{|y|\lesssim 2^{-j+M}}\Big(\cdots\Big)dy+\sum_{l\geq 0}\int_{|y|\sim 2^{-j+M+l}}\Big(\cdots\Big)dy\\
&\lesssim\int_{|y|\lesssim 2^{-j+M}}\prod_{i=1}^{n}\Big(2^{Md(\frac{1}{q_i}-\frac{1}{2})}\mathrm{M}_{q_i}(f_i*\Phi_j^i)(x)\Big)\omega_k^N(y)dy\\
&\quad+\sum_{l\geq 0}\prod_{i=1}^{n}\Big(2^{\frac{ld}{2}}2^{(M+l)d(\frac{1}{q_i}-\frac{1}{2})}\mathrm{M}_{q_i}(f_i*\Phi_j^i)(x)\Big)\Big(\frac{1}{(1+2^{k-j+M+l})^{N-d-1}}\Big)\\
&\lesssim\prod_{i=1}^{n}\Big(2^{Md(\frac{1}{q_i}-\frac{1}{2})}\mathrm{M}_{q_i}(f_i*\Phi_j^i)(x)\Big),
\end{aligned}\tag{6.6}$$

where we use $\omega_k^N(y)=2^{kd}(1+|2^ky|)^{-N}$ for the second inequality.

By (6.3) and (6.6)

$$\begin{aligned}
|\mathrm{A}_{\mathbf{m}_\lambda}^{M,j,k}(f_1,\cdots,f_n)(x)|&\lesssim\frac{2^{-Ms+\frac{Mnd}{2}}2^{j\delta}}{2^k}\sup_{x\in\mathbb{R}^d}\Big(\sum_{|\alpha|\leq 1}2^{-j\delta}\big\|\partial_x^\alpha\mathbf{m}_\lambda(x,2^j\vec{\cdot})\widehat{\Psi}(\vec{\cdot})\big\|_{L_s^2((\mathbb{R}^d)^n)}\Big)\\
&\quad\times\prod_{i=1}^{n}\Big(2^{Md(\frac{1}{q_i}-\frac{1}{2})}\mathrm{M}_{q_i}(f_i*\Phi_j^i)(x)\Big).
\end{aligned}\tag{6.7}$$

Since $\frac{p}{p_1}+\cdots+\frac{p}{p_n}=1$, by Hölder's inequality, if $0<q_i<\min(2,p_i)$, then we have

$$\Big\|\prod_{i=1}^{n}\mathrm{M}_{q_i}(f_i*\Phi_j^i)\Big\|_{L^p}\lesssim\prod_{i=1}^{n}\|\mathrm{M}_{q_i}(f_i*\Phi_j^i)\|_{L^{p_i}}\lesssim\prod_{i=1}^{n}\|f_i*\Phi_j^i\|_{L^{p_i}}\lesssim\prod_{i=1}^{n}\|f_i\|_{h^{p_i}}.\tag{6.8}$$

Since

$$\|\mathrm{I}_{\mathbf{m}_\lambda}^{N_1,N_2}(f_1,\cdots,f_n)\|_{L^p}^{\min(1,p)}\leq\sum_{j=0}^{\infty}\sum_{k=j+10}^{\infty}\sum_{M\geq 0}\|\mathrm{A}_{\mathbf{m}_\lambda}^{M,j,k}(f_1,\cdots,f_n)\|_{L^p}^{\min(1,p)},$$

by (6.7) and (6.8), if $0<q_i<\min(2,p_i)$ for all $1\leq i\leq n$, then

$$\begin{aligned}
\|\mathrm{I}_{\mathbf{m}_\lambda}^{N_1,N_2}(f_1,\cdots,f_n)\|_{L^p}&\lesssim\sup_{j\geq 0}\sup_{x\in\mathbb{R}^d}\Big(\sum_{|\alpha|\leq 1}2^{-j\delta}\big\|\partial_x^\alpha\mathbf{m}_\lambda(x,2^j\vec{\cdot})\widehat{\Psi}(\vec{\cdot})\big\|_{L_s^2((\mathbb{R}^d)^n)}\Big)\prod_{i=1}^{n}\|f_i\|_{h^{p_i}}\\
&\lesssim\sup_{j\geq 0}\sup_{x\in\mathbb{R}^d}\Big(\sum_{|\alpha|\leq 1}2^{-j\delta}\big\|\partial_x^\alpha\mathbf{m}(x,2^j\vec{\cdot})\widehat{\Psi}(\vec{\cdot})\big\|_{L_s^2((\mathbb{R}^d)^n)}\Big)\prod_{i=1}^{n}\|f_i\|_{h^{p_i}},
\end{aligned}\tag{6.9}$$

uniformly in $0<\lambda\leq 1$, $N_1$, and $N_2$ when $\frac{s}{d}>\sum_{i=1}^{n}\frac{1}{q_i}$. Therefore, by taking $q_i\nearrow\min(2,p_i)$ for $1\leq i\leq n$, we have (6.9) when

$$\Big(\frac{1}{p_1},\cdots,\frac{1}{p_n}\Big)\in\mathrm{B}_n\Big(\frac{s}{d}\Big)=\Big\{(x_1,\cdots,x_n)\in(0,\infty)^n:\sum_{i=1}^{n}\max(x_i,1/2)<\frac{s}{d}\Big\}.$$

## 7. Proof of Theorem 1.2 : Estimates for the term $\mathrm{II}^{N_1}_{\mathbf{m}_\lambda}(f_1,\cdots,f_n)$

Recall the definition

$$\mathscr{T}^j_{\mathbf{m}_\lambda}(f_1,\cdots,f_n)(x)=\gamma(\lambda x)\int_{(\mathbb{R}^d)^n}\mathbf{m}(x,\vec{\xi})\widehat{\Psi}(2^{-j}\vec{\xi})\Big(\prod_{i=1}^n\big(\widehat{\Phi}^i_j(\xi_i)\widehat{f_i}(\xi_i)\big)\Big)\widehat{\Phi}^{n+1}_j(\xi_{n+1})e^{2\pi i\sum_{i=1}^n\langle x,\xi_i\rangle}d\vec{\xi}.$$

To estimate the term $\mathrm{II}^{N_1}_{\mathbf{m}_\lambda}(f_1,\cdots,f_n)$ in (4.3), we consider it in two cases:

(1) $\widehat{\Phi}^{n+1}$ is compactly supported away from the origin,
(2) $\widehat{\Phi}^{n+1}$ is not compactly supported away from the origin.

### 7.1. Proof of Theorem 1.2 : Estimates for the term $\mathrm{II}^{N_1}_{\mathbf{m}_\lambda}(f_1,\cdots,f_n)$ in (4.3) when $\widehat{\Phi}^{n+1}$ is not compactly supported away from the origin.

If $\widehat{\Phi}^{n+1}$ is not compactly supported away from the origin. Then by Lemma 3.11, there are two indices $i_1,i_2\in\{1,2,\cdots,n\}$ so that $\widehat{\Phi}^{i_1}$ and $\widehat{\Phi}^{i_2}$ are compactly supported away from the origin. Without loss of generality let $i_1=1$ and $i_2=2$. Recall the definition

$$\mathrm{II}^{N_1}_{\mathbf{m}_\lambda}(f_1,\cdots,f_n)(x)=\sum_{j=0}^{N_1}\mathscr{T}^j_{\mathbf{m}_\lambda}(f_1,\cdots,f_n)*\varphi_{j-10}(x).$$

For notational convenience, we use $\varphi_j$ instead of $\varphi_{j-10}$. Then if we follow the estimates for $\mathscr{T}^j_{\mathbf{m}_\lambda}(f_1,\cdots,f_n)*\psi_j(x)$ in Section 5 we have

$$\begin{aligned}&|\mathscr{T}^j_{\mathbf{m}_\lambda}(f_1,\cdots,f_n)*\varphi_j(x)|\\&\lesssim\Big(\sum_{M\geq 0}2^{-Ms+\frac{Mnd}{2}+\sum_{i=1}^n Md\left(\frac{1}{q_i}-\frac12\right)}\Big)\Big(\sup_{x\in\mathbb{R}^d}\big\|\mathbf{m}(x,2^j\vec{\cdot})\widehat{\Psi}(\vec{\cdot})\big\|_{L^2_s((\mathbb{R}^d)^n)}\Big)\Big(\prod_{i=1}^n\mathrm{M}_{q_i}(f_i*\Phi^i_j)(x)\Big).\end{aligned}$$

Since $\frac{p}{p_1}+\cdots+\frac{p}{p_n}=1$, by Hölder's inequality, if $0<q_i<\min(2,p_i)$, then by Lemma 3.4 we have

(7.1)

$$\begin{aligned}&\Big\|\mathrm{II}^{N_1}_{\mathbf{m}_\lambda}(f_1,\cdots,f_n)\Big\|_{L^p}\\&\lesssim\Big(\sum_{M\geq0}2^{-Ms+\frac{Mnd}{2}+\sum_{i=1}^n Md\left(\frac{1}{q_i}-\frac12\right)}\Big)\Big(\sup_{j\geq0}\sup_{x\in\mathbb{R}^d}\big\|\mathbf{m}(x,2^j\vec{\cdot})\widehat{\Psi}(\vec{\cdot})\big\|_{L^2_s((\mathbb{R}^d)^n)}\Big)\\&\times\Big\|\Big(\sum_{j=0}^{N_1}|\mathrm{M}_{q_1}(f_1*\Phi^1_j)|^2\Big)^{\frac12}\Big\|_{L^{p_1}}\Big\|\Big(\sum_{j=0}^{N_1}|\mathrm{M}_{q_2}(f_2*\Phi^2_j)|^2\Big)^{\frac12}\Big\|_{L^{p_2}}\Big(\prod_{i=3}^n\big\|\mathrm{M}_{q_i}(\sup_{j\geq0}|f_i*\Phi^i_j|)\big\|_{L^{p_i}}\Big)\\&\lesssim\Big(\sup_{j\geq0}\sup_{x\in\mathbb{R}^d}\big\|\mathbf{m}(x,2^j\vec{\cdot})\widehat{\Psi}(\vec{\cdot})\big\|_{L^2_s((\mathbb{R}^d)^n)}\Big)\Big\|\Big(\sum_{j=0}^{N_1}|f_1*\Phi^1_j|^2\Big)^{\frac12}\Big\|_{L^{p_1}}\Big\|\Big(\sum_{j=0}^{N_1}|f_2*\Phi^2_j|^2\Big)^{\frac12}\Big\|_{L^{p_2}}\Big(\prod_{i=3}^n\big\|\sup_{j\geq0}|f_i*\Phi^i_j|\big\|_{L^{p_i}}\Big)\\&\lesssim\Big(\sup_{j\geq0}\sup_{x\in\mathbb{R}^d}\big\|\mathbf{m}(x,2^j\vec{\cdot})\widehat{\Psi}(\vec{\cdot})\big\|_{L^2_s((\mathbb{R}^d)^n)}\Big)\prod_{i=1}^n\|f_i\|_{h^{p_i}}\end{aligned}$$

uniformly in $0<\lambda\leq1$ and $N_1$ when $\frac{s}{d}>\sum_{i=1}^n\frac{1}{q_i}$. Therefore, by taking $q_i\nearrow\min(2,p_i)$ for $1\leq i\leq n$, we have (7.1) when

$$\Big(\frac{1}{p_1},\cdots,\frac{1}{p_n}\Big)\in\mathrm{B}_n\Big(\frac{s}{d}\Big)=\Big\{(x_1,\cdots,x_n)\in(0,\infty)^n:\sum_{i=1}^n\max(x_i,1/2)<\frac{s}{d}\Big\}.$$

7.2. **Proof of Theorem 1.2 : Estimates for the term $II^{N_1}_{\mathbf{m}_\lambda}(f_1,\cdots,f_n)$ in (4.3) when $\widehat{\Phi}^{n+1}$ is compactly supported away from the origin.** In this case the estimates for the term $II^{N_1}_{\mathbf{m}_\lambda}(f_1,\cdots,f_n)$ are similar to those for $I^{N_1,N_2}_{\mathbf{m}_\lambda}(f_1,\cdots,f_n)$. Recall the definition

$$II^{N_1}_{\mathbf{m}_\lambda}(f_1,\cdots,f_n)(x)=\sum_{j=0}^{N_1}\mathscr{T}^j_{\mathbf{m}_\lambda}(f_1,\cdots,f_n)*\varphi_{j-10}(x).$$

For notational convenience, we use $\varphi_j$ instead of $\varphi_{j-10}$. Note that

$$\begin{aligned}&\mathscr{T}^j_{\mathbf{m}_\lambda}(f_1,\cdots,f_n)*\varphi_j(x)\\&=\int_{\mathbb{R}^{(n+2)d}}\mathbf{m}_\lambda(x-y,\vec{\xi})\widehat{\Psi}(2^{-j}\vec{\xi})\,\widehat{\Phi}^{n+1}_j(-\xi_1-\cdots-\xi_n)\\&\quad\times\Big(\prod_{i=1}^n\big(\widehat{\Phi}^i_j(\xi_i)\widehat{f}_i(\xi_i)\big)\Big)\widehat{\varphi}_j(\eta)e^{2\pi i\langle y,\eta-(\xi_1+\cdots+\xi_n)\rangle}e^{2\pi i\langle x,\xi_1+\cdots+\xi_n\rangle}dyd\eta d\vec{\xi}.\end{aligned}$$

Since $\widehat{\Phi}^{n+1}$ is compactly supported away from the origin, $\widehat{\Phi}^{n+1}_j(-\xi_1-\cdots-\xi_n)\,\widehat{\varphi}_j(\eta)\neq 0$ only if $|\xi_1+\cdots+\xi_n-\eta|\sim 2^j$. And the term $\widehat{\Phi}^{n+1}_j(-\xi_1-\cdots-\xi_n)\,\widehat{\varphi}_j(\eta)$ can be written as a finite sum of the form:

$$\widehat{\Phi}^{n+1}_j(-\xi_1-\cdots-\xi_n)\,\widehat{\varphi}_j(\eta)\prod_{i=1}^d\phi^i\Big(\frac{\xi_1^i+\cdots+\xi_n^i-\eta^i}{2^j}\Big)$$

where $\phi^i(t)$ are smooth functions that are supported in $|t|\lesssim 1$ and at least one of $\phi^i(t)$ is supported in $|t|\sim 1$. Thus $\mathscr{T}^j_{\mathbf{m}_\lambda}(f_1,\cdots,f_n)*\varphi_j(x)$ can be written as a finite sum of the form:

$$\begin{aligned}&\mathscr{T}^j_{\mathbf{m}_\lambda}(f_1,\cdots,f_n)*\varphi_j(x)\\&:=\int_{\mathbb{R}^{(n+2)d}}\prod_{i=1}^d\phi^i\Big(\frac{\xi_1^i+\cdots+\xi_n^i-\eta^i}{2^j}\Big)\mathbf{m}_\lambda(x-y,\vec{\xi})\widehat{\Psi}(2^{-j}\vec{\xi})\,\widehat{\Phi}^{n+1}_j(-\xi_1-\cdots-\xi_n)\\&\quad\times\Big(\prod_{i=1}^n\big(\widehat{\Phi}^i_j(\xi_i)\widehat{f}_i(\xi_i)\big)\Big)\widehat{\varphi}_j(\eta)e^{2\pi i\langle y,\eta-(\xi_1+\cdots+\xi_n)\rangle}e^{2\pi i\langle x,\xi_1+\cdots+\xi_n\rangle}dyd\eta d\vec{\xi}.\end{aligned}$$

where $\phi^l(t)$ is supported in $|t|\sim 1$ for some $1\le l\le d$. Then by using integration by parts via

$$\partial_{y_l}e^{2\pi i\langle y,\,\eta-(\xi_1+\cdots+\xi_n)\rangle}=2\pi i(\eta^l-\xi_1^l-\cdots-\xi_n^l)e^{2\pi i\langle y,\,\eta-(\xi_1+\cdots+\xi_n)\rangle},$$

we obtain that

$$\begin{aligned}&\mathscr{T}^j_{\mathbf{m}_\lambda}(f_1,\cdots,f_n)*\varphi_j(x)\\&=-\int_{\mathbb{R}^{(n+2)d}}\frac{e^{2\pi i\langle x,\,\xi_1+\cdots+\xi_n\rangle}}{2\pi i(\eta^l-\xi_1^l-\cdots-\xi_n^l)}\prod_{i=1}^d\phi^i\Big(\frac{\xi_1^i+\cdots+\xi_n^i-\eta^i}{2^j}\Big)\big[\partial_{y_l}\mathbf{m}_\lambda(x-y,\vec{\xi})\widehat{\Psi}(2^{-j}\vec{\xi})\big]\\&\quad\times\Big(\prod_{i=1}^n\big(\widehat{\Phi}^i_j(\xi_i)\widehat{f}_i(\xi_i)\big)\Big)\widehat{\Phi}^{n+1}_j(-\xi_1-\cdots-\xi_n)\widehat{\varphi}_j(\eta)e^{2\pi i\langle y,\,\eta-(\xi_1+\cdots+\xi_n)\rangle}dyd\eta d\vec{\xi}.\end{aligned}$$

Then by using the identity

$$\partial_{y_l}\mathbf{m}_\lambda(x-y,\vec{\xi})\widehat{\Psi}(2^{-j}\vec{\xi})=\int_{(\mathbb{R}^d)^n}\mathscr{F}^{-1}\Big(\partial_{y_l}\big[\mathbf{m}_\lambda(x-y,2^j\vec{\cdot})\widehat{\Psi}(\vec{\cdot})\big]\Big)(\vec{z})e^{-2\pi i\langle\vec{z},2^{-j}\vec{\xi}\rangle}d\vec{z},$$

we have

$$\mathscr{T}^j_{\mathbf{m}_\lambda}(f_1,\cdots,f_n)*\varphi_j(x)$$
$$(7.2)\quad = -\int_{\mathbb{R}^{(2n+2)d}} \frac{e^{2\pi i\langle(\vec{\xi},\eta),(x-y-2^{-j}z_1,\cdots,x-y-2^{-j}z_n,y)\rangle}}{2\pi i(\eta^l-\xi_1^l-\cdots-\xi_n^l)} \prod_{i=1}^d \phi^i\Big(\frac{\xi_1^i+\cdots+\xi_n^i-\eta^i}{2^j}\Big)\widehat{\Phi}^{n+1}_j(-\xi_1-\cdots-\xi_n)$$
$$\times \mathscr{F}^{-1}\Big(\partial_{y_l}\big[\mathbf{m}_\lambda(x-y,2^j\vec{\cdot})\widehat{\Psi}(\vec{\cdot})\big]\Big)(\vec{\mathbf{z}})\Big(\prod_{i=1}^n\big(\widehat{\Phi}^i_j(\xi_i)\widehat{f}_i(\xi_i)\big)\Big)\widehat{\varphi}_j(\eta)dyd\eta d\vec{\mathbf{z}}d\vec{\xi}.$$

By (7.2) and Lemma 6.1, we have

$$|\mathscr{T}^j_{\mathbf{m}_\lambda}(f_1,\cdots,f_n)*\varphi_j(x)| \le C_N\frac{1}{2^j}\iint \Big|\mathscr{F}^{-1}\Big(\partial_{y_l}\big[\mathbf{m}_\lambda(x-y,2^j\vec{\cdot})\widehat{\Psi}(\vec{\cdot})\big]\Big)(\vec{\mathbf{z}})\Big|$$
$$\times\Big(\prod_{i=1}^n \omega_j^N * |f_i*\Phi^i_j|(x-y-2^{-j}z_i)\Big)(\omega_j^N(y))dyd\vec{\mathbf{z}}.$$

Let

$$\Gamma_0 := \big\{\vec{\mathbf{z}}\in(\mathbb{R}^d)^n : |\vec{\mathbf{z}}|\le 1\big\},\quad \Gamma_M := \big\{\vec{\mathbf{z}}\in(\mathbb{R}^d)^n : 2^{M-1}<|\vec{\mathbf{z}}|\le 2^M\big\},\quad M\ge 1.$$

Then we have

$$|\mathscr{T}^j_{\mathbf{m}_\lambda}(f_1,\cdots,f_n)*\varphi_j(x)| \lesssim \sum_{M\ge 0} \mathrm{B}^{M,j}_{\mathbf{m}_\lambda}(f_1,\cdots,f_n)(x)$$

where

$$(7.3)\quad \mathrm{B}^{M,j}_{\mathbf{m}_\lambda}(f_1,\cdots,f_n)(x) := \frac{1}{2^j}\int_{\mathbb{R}^d}\int_{\Gamma_M} \Big|\mathscr{F}^{-1}\Big(\partial_{y_l}\big[\mathbf{m}_\lambda(x-y,2^j\vec{\cdot})\widehat{\Psi}(\vec{\cdot})\big]\Big)(\vec{\mathbf{z}})\Big|$$
$$\times\Big(\prod_{i=1}^n \omega_j^N * |f_i*\Phi^i_j|(x-y-2^{-j}z_i)\Big)(\omega_j^N(y))d\vec{\mathbf{z}}dy.$$

Since the case $M=0$ is similar to the case $M\ge 1$, we only consider the case $M\ge 1$. By (7.3) we have

$$\mathrm{B}^{M,j}_{\mathbf{m}_\lambda}(f_1,\cdots,f_n)(x) \lesssim \frac{2^{Mnd}}{2^j}\int_{\mathbb{R}^d}\int_{|\vec{\mathbf{z}}|\sim 1} \Big|\mathscr{F}^{-1}\Big(\partial_{y_l}\big[\mathbf{m}_\lambda(x-y,2^j\vec{\cdot})\widehat{\Psi}(\vec{\cdot})\big]\Big)(2^M\vec{\mathbf{z}})\Big|$$
$$\times\Big(\prod_{i=1}^n \omega_j^N * |f_i*\Phi^i_j|(x-y-2^{-j+M}z_i)\Big)(\omega_j^N(y))d\vec{\mathbf{z}}dy.$$

By using Hölder's inequality with the $\vec{\mathbf{z}}$ variable

$$(7.4)\quad |\mathrm{B}^{M,j}_{\mathbf{m}_\lambda}(f_1,\cdots,f_n)(x)| \lesssim \frac{2^{-Ms+\frac{Mnd}{2}}2^{j\delta}}{2^j}\sup_{x\in\mathbb{R}^d}\Big(\sum_{|\alpha|\le 1} 2^{-j\delta}\big\|\partial_x^\alpha\mathbf{m}_\lambda(x,2^j\vec{\cdot})\widehat{\Psi}(\vec{\cdot})\big\|_{L^2_s((\mathbb{R}^d)^n)}\Big)$$
$$\times\int\Big(\prod_{i=1}^n\Big(\int_{|z_i|\lesssim 1}\big|\omega_j^N*|f_i*\Phi^i_j|(x-y-2^{-j+M}z_i)\big|^2dz_i\Big)^{1/2}\Big)\omega_j^N(y)dy.$$

Let $|y|\sim 2^{-j+M+l}$ for some $l\geq 0$, then $y=2^{-j+M+l}y'$ for some $|y'|\lesssim 1$. And by the change of variables $y'+2^{-l}z_i\to z_i$ we have

$$\begin{aligned}&\Big(\int_{|z_i|\lesssim 1}\big|\omega_j^N*|f_i*\Phi_j^i|(x-y-2^{-j+M}z_i)\big|^2dz_i\Big)^{1/2}\\ &=\Big(\int_{|z_i|\lesssim 1}\big|\omega_j^N*|f_i*\Phi_j^i|(x-2^{-j+M+l}(y'+2^{-l}z_i))\big|^2dz_i\Big)^{1/2}\\ &=2^{\frac{ld}{2}}\Big(\int_{|z_i|\lesssim 1}\big|\omega_j^N*|f_i*\Phi_j^i|(x-2^{-j+M+l}z_i)\big|^2dz_i\Big)^{1/2}.\end{aligned}\tag{7.5}$$

By Lemma 3.9, for any $0<q_i\leq 2$, if $N>(\frac{2}{q_i}+1)d$, then

$$\Big(\int_{|z_i|\lesssim 1}\big|\omega_j^N*|f_i*\Phi_j^i|(x-2^{-j+M}z_i)\big|^2dz_i\Big)^{1/2}\lesssim 2^{(M+l)d(\frac{1}{q_i}-\frac12)}\mathrm{M}_{q_i}(f_i*\Phi_j^i)(x).\tag{7.6}$$

Therefore by (7.5) and (7.6), if we take $N$ large enough, then

$$\begin{aligned}&\int\Big(\prod_{i=1}^n\Big(\int_{|z_i|\lesssim 1}\big|\omega_j^N*|f_i*\Phi_j^i|(x-y-2^{-j+M}z_i)\big|^2dz_i\Big)^{1/2}\Big)\omega_j^N(y)dy\\ &\lesssim\int_{|y|\lesssim 2^{-j+M}}\Big(\cdots\Big)dy+\sum_{l\geq 0}\int_{|y|\sim 2^{-j+M+l}}\Big(\cdots\Big)dy\\ &\lesssim\int_{|y|\lesssim 2^{-j+M}}\Big(\prod_{i=1}^n 2^{Md(\frac{1}{q_i}-\frac12)}\mathrm{M}_{q_i}(f_i*\Phi_j^i)(x)\Big)\omega_j^N(y)dy\\ &\quad+\sum_{l\geq 0}\prod_{i=1}^n\Big(2^{\frac{ld}{2}}2^{(M+l)d(\frac{1}{q_i}-\frac12)}\mathrm{M}_{q_i}(f_i*\Phi_j^i)(x)\Big)\Big(\frac{1}{(1+2^{M+l})^{N-d-1}}\Big)\\ &\lesssim\prod_{i=1}^n\Big(2^{Md(\frac{1}{q_i}-\frac12)}\mathrm{M}_{q_i}(f_i*\Phi_j^i)(x)\Big)\end{aligned}\tag{7.7}$$

where we use $\omega_j^N(y)=2^{jd}(1+|2^jy|)^{-N}$ for the second inequality.

By (7.4) and (7.7)

$$\begin{aligned}|\mathrm{B}^{M,j}_{\mathbf{m}_\lambda}(f_1,\cdots,f_n)(x)|\lesssim&\frac{2^{-Ms+\frac{Mnd}{2}}2^{j\delta}}{2^j}\sup_{x\in\mathbb{R}^d}\Big(\sum_{|\alpha|\leq 1}2^{-j\delta}\big\|\partial_x^\alpha\mathbf{m}_\lambda(x,2^j\,\vec{\cdot}\,)\widehat{\Psi}(\vec{\cdot}\,)\big\|_{L^2_s((\mathbb{R}^d)^n)}\Big)\\ &\times\prod_{i=1}^n\Big(2^{Md(\frac{1}{q_i}-\frac12)}\mathrm{M}_{q_i}(f_i*\Phi_j^i)(x)\Big).\end{aligned}\tag{7.8}$$

Since $\frac{p}{p_1}+\cdots+\frac{p}{p_n}=1$, by Hölder's inequality, if $0<q_i<\min(2,p_i)$, then we have

$$\Big\|\prod_{i=1}^n\mathrm{M}_{q_i}(f_i*\Phi_j^i)\Big\|_{L^p}\lesssim\prod_{i=1}^n\|\mathrm{M}_{q_i}(f_i*\Phi_j^i)\|_{L^{p_i}}\lesssim\prod_{i=1}^n\|f_i*\Phi_j^i\|_{L^{p_i}}\lesssim\prod_{i=1}^n\|f_i\|_{h^{p_i}}.\tag{7.9}$$

Since

$$\|\mathrm{II}^{N_1}_{\mathbf{m}_\lambda}(f_1,\cdots,f_n)\|_{L^p}^{\min(1,p)}\leq\sum_{j=0}^\infty\sum_{M\geq 0}\|\mathrm{B}^{M,j}_{\mathbf{m}_\lambda}(f_1,\cdots,f_n)\|_{L^p}^{\min(1,p)},$$

by (7.8) and (7.9), if $0 < q_i < \min(2, p_i)$ for all $1 \le i \le n$, then

$$\begin{aligned}\|\mathrm{II}^{N_1}_{\mathbf{m}_\lambda}(f_1,\cdots,f_n)\|_{L^p} &\lesssim \sup_{j\ge 0}\sup_{x\in\mathbb{R}^d}\Big(\sum_{|\alpha|\le 1} 2^{-j\delta}\big\|\partial_x^\alpha \mathbf{m}_\lambda(x,2^j\vec{\cdot})\widehat{\Psi}(\vec{\cdot})\big\|_{L^2_s((\mathbb{R}^d)^n)}\Big)\prod_{i=1}^n \|f_i\|_{h^{p_i}} \\ &\lesssim \sup_{j\ge 0}\sup_{x\in\mathbb{R}^d}\Big(\sum_{|\alpha|\le 1} 2^{-j\delta}\big\|\partial_x^\alpha \mathbf{m}(x,2^j\vec{\cdot})\widehat{\Psi}(\vec{\cdot})\big\|_{L^2_s((\mathbb{R}^d)^n)}\Big)\prod_{i=1}^n \|f_i\|_{h^{p_i}},\end{aligned} \tag{7.10}$$

uniformly in $0 < \lambda \le 1$ and $N_1$ when $\frac{s}{d} > \sum_{i=1}^n \frac{1}{q_i}$. Therefore, by taking $q_i \nearrow \min(2, p_i)$ for $1 \le i \le n$, we have (7.10) when

$$\Big(\frac{1}{p_1},\cdots,\frac{1}{p_n}\Big) \in \mathrm{B}_n\Big(\frac{s}{d}\Big) = \Big\{(x_1,\cdots,x_n)\in(0,\infty)^n : \sum_{i=1}^n \max(x_i,1/2) < \frac{s}{d}\Big\}.$$

## 8. Proof of Theorem 1.2 : Estimates for the term $\mathscr{T}^{-1}_{\mathbf{m}_\lambda}(f_1,\cdots,f_n)$ in (4.1)

Recall the definition

$$\mathscr{T}^{-1}_{\mathbf{m}_\lambda}(f_1,\cdots,f_n)(x) := \gamma(\lambda x)\int_{(\mathbb{R}^d)^n} \mathbf{m}(x,\vec{\xi})\widehat{\Phi}(\vec{\xi})\Big(\prod_{i=1}^n \widehat{f_i}(\xi_i)\Big)e^{2\pi i\sum_{i=1}^n\langle x,\xi_i\rangle}d\vec{\xi}.$$

Let $\varphi$ and $\psi$ be as in (3.13), then

$$\sum_{k=10}^{\infty}\widehat{\psi_k}(\eta) + \widehat{\varphi}(2^{-10}\eta) = 1.$$

Thus we decompose $\mathscr{T}^{-1}_{\mathbf{m}_\lambda}(f_1,\cdots,f_n)(x)$ into

$$\mathrm{IV}_{\mathbf{m}_\lambda}(f_1,\cdots,f_n)(x) + \mathrm{V}_{\mathbf{m}_\lambda}(f_1,\cdots,f_n)(x)$$

for almost every $x$, where $\mathrm{IV}_{\mathbf{m}_\lambda}(f_1,\cdots,f_n)$ and $\mathrm{V}_{\mathbf{m}_\lambda}(f_1,\cdots,f_n)$ are $L^2(\mathbb{R}^d)$ functions given by

$$\begin{aligned}\mathrm{IV}_{\mathbf{m}_\lambda}(f_1,\cdots,f_n)(x) &:= \mathscr{F}^{-1}\Big(\sum_{k=10}^{\infty}\mathscr{F}\big(\mathscr{T}^{-1}_{\mathbf{m}_\lambda}(f_1,\cdots,f_n)\big)(\eta)\,\widehat{\psi}(2^{-k}\eta)\Big)(x), \\ \mathrm{V}_{\mathbf{m}_\lambda}(f_1,\cdots,f_n)(x) &:= \mathscr{F}^{-1}\Big(\mathscr{F}\big(\mathscr{T}^{-1}_{\mathbf{m}_\lambda}(f_1,\cdots,f_n)\big)(\eta)\,\widehat{\varphi}(2^{-10}\eta)\Big)(x).\end{aligned} \tag{8.1}$$

The estimates for the term $\mathrm{IV}_{\mathbf{m}_\lambda}(f_1,\cdots,f_n)$ are similar to those for the term $\mathrm{I}^{0,N_2}_{\mathbf{m}_\lambda}(f_1,\cdots,f_n)$ in (4.3). And the estimates for the term $\mathrm{V}_{\mathbf{m}_\lambda}(f_1,\cdots,f_n)$ are similar to those for the term $\mathrm{II}^{0}_{\mathbf{m}_\lambda}(f_1,\cdots,f_n)$ in (4.3). Detailed estimates will be omitted.

## 9. Appendix

**Proof of Lemma 3.4.** For the proof of Lemma 3.4, we adopt the proof of Theorem 2.2.9 in [14]. Let $\Phi \in \mathscr{S}(\mathbb{R}^d)$ with $\int \Phi = 1$. Let $f \in h^p \cap L^1$ and $M \in \mathbb{Z}^+$. Let $r_j$ be the Rademacher functions(see Appendix C.1 in [13]). Note that

$$\Big|\sum_{j=0}^M r_j(w)\Psi_j * f\Big| \le \sup_{0<\epsilon<1}\Big|\Phi_\epsilon * \Big(\sum_{j=0}^M r_j(w)\Psi_j * f\Big)\Big|$$

which holds since $\{\Phi_\epsilon\}_{0<\epsilon<1}$ is an approximate identity. We this inequality to the power $p$, and integrate over $x \in \mathbb{R}^d$ and $w \in [0,1]$, and we use the maximal characterization of $h^p$ to obtain

$$\int_0^1\int_{\mathbb{R}^d}\Big|\sum_{j=0}^M r_j(w)\Psi_j * f(x)\Big|^p dx\,dw \le C^p_{p,d}\int_0^1\Big\|\sum_{j=0}^M r_j(w)\Psi_j * f\Big\|^p_{h^p} dw. \tag{9.1}$$

Let $K(x) := \sum_{j=0}^{M} r_j(w)\Psi_j(x)$. Then from

$$\widehat{K}(\xi) = \sum_{j=0}^{M} r_j(w)\widehat{\Psi}(\xi/2^j)$$

we see that $\widehat{K}$ belongs to $\mathscr{S}^0_{1,0}(\mathbb{R}^d)$ unformly in $w \in [0,1]$. Thus by Theorem B we obtain

$$\Big\| \sum_{j=0}^{M} r_j(w)\Psi_j * f \Big\|_{L^p}^p \lesssim \Big\| \sum_{j=0}^{M} r_j(w)\Psi_j * f \Big\|_{h^p}^p \lesssim \|f\|_{h^p}. \tag{9.2}$$

And by Khintchine's inequality(see Appendix C.2 in [14]) we have

$$\int_{\mathbb{R}^d} \Big( \sum_{j=0}^{M} |\psi_j * f|^2 \Big)^{\frac{p}{2}} dx \lesssim \int_0^1 \int_{\mathbb{R}^d} \Big| \sum_{j=0}^{M} r_j(w)\Psi_j * f(x) \Big|^p dx\, dw. \tag{9.3}$$

By (9.1), (9.2), and (9.3) we have

$$\Big\| \Big( \sum_{j=0}^{M} |f * \Psi_j|^2 \Big)^{\frac{1}{2}} \Big\|_{L^p(\mathbb{R}^d)} \le C_{d,p,\Psi} \|f\|_{h^p(\mathbb{R}^d)}.$$

By letting $M \to \infty$ we have the desired results for $f \in h^p \cap L^1$. Since $h^p \cap C^\infty$ is dense in $h^p$(see [11, 12]). Using density, we can extend this estimate to all $f \in h^p$.

**Proof of Lemma 1.1.** We first observe that the condition (2) in Theorem A is equivalent to

$$\sum_{i \in I^c} \frac{1}{p_i} < \Big(\frac{s}{d} + \frac{1}{2}\Big) - \frac{\sharp I}{2} = \Big(\frac{s}{d} + \frac{1}{2}\Big) - \frac{n - \sharp I^c}{2}. \tag{9.4}$$

By replacing the set $I^c = J_n \setminus I$ in (9.4) with $I$, the collection of $n$-tuples $(1/p_1, \cdots, 1/p_n)$ for which the condition (2) holds is equivalent to the set $\mathrm{B}_n\big(\frac{s}{d} + \frac{1}{2}\big)$ where

$$\mathrm{B}_n(\alpha) := \bigcap_{I \subset J_n} \Big\{ (x_1, \cdots, x_n) \in (0,\infty)^n \,:\, \sum_{i \in I} \Big(x_i - \frac{1}{2}\Big) + \frac{n}{2} < \alpha \Big\}.$$

Then we claim that $\mathrm{B}_n(\alpha) = \mathrm{A}_n(\alpha)$ for $\alpha > 0$ where

$$\mathrm{A}_n(\alpha) := \Big\{ (x_1, \cdots, x_n) \in (0,\infty)^n \,:\, \sum_{i=1}^{n} \max\Big(x_i, \frac{1}{2}\Big) < \alpha \Big\}.$$

To see this, for each $I \subset J_n$ we set

$$R_I := \Big\{ (x_1, \cdots, x_n) \in (0,\infty)^n : x_i > 1/2 \quad \text{if } i \in I \ \text{ and } \ x_i \le 1/2 \ \text{ if } i \in I^c \Big\}.$$

Then we have $(0,\infty)^n = \bigcup_{I' \subset J_n} R_{I'}$. First we prove that $\mathrm{A}_n(\alpha) \subset \mathrm{B}_n(\alpha)$. Let $(x_1, \cdots, x_n) \in \mathrm{A}_n(\alpha) \cap R_{I'}$, then

$$\sum_{i=1}^{n} \max\Big(x_i, \frac{1}{2}\Big) = \sum_{i \in I'} \max\Big(x_i, \frac{1}{2}\Big) + \sum_{i \in (I')^c} \max\Big(x_i, \frac{1}{2}\Big) = \sum_{i \in I'} x_i + \frac{n - \sharp I'}{2} < \alpha. \tag{9.5}$$

Let $I \subset J_n$, then since $\sum_{i \in I \cap (I')^c} (x_i - \frac{1}{2}) \le 0$, by (9.5) we have

$$\begin{aligned} \sum_{i \in I} \Big(x_i - \frac{1}{2}\Big) + \frac{n}{2} &= \sum_{i \in I \cap I'} \Big(x_i - \frac{1}{2}\Big) + \sum_{i \in I \cap (I')^c} \Big(x_i - \frac{1}{2}\Big) + \frac{n}{2} \\ &\le \sum_{i \in I \cap I'} \Big(x_i - \frac{1}{2}\Big) + \frac{n}{2} \le \sum_{i \in I'} \Big(x_i - \frac{1}{2}\Big) + \frac{n}{2} < \alpha, \end{aligned}$$

which implies that $\mathrm{A}_n(\alpha) \cap R_{I'} \subset \mathrm{B}_n(\alpha)$, and so $\mathrm{A}_n(\alpha) \subset \mathrm{B}_n(\alpha)$.

Conversely, let $(x_1, \cdots, x_n) \in \mathrm{B}_n(\alpha) \cap R_{I'}$, then for each $I \subset J_n$ we have

$$\sum_{i \in I} \Big(x_i - \frac{1}{2}\Big) + \frac{n}{2} < \alpha.$$

Then by taking $I = I'$, we have

$$\sum_{i \in I'} \Big(x_i - \frac{1}{2}\Big) + \frac{n}{2} < \alpha,$$

and by (9.5) this implies that $\mathrm{B}_n(\alpha) \cap R_{I'} \subset \mathrm{A}_n(\alpha) \cap R_{I'}$. □

YARYONG HEO, DEPARTMENT OF MATHEMATICS, KOREA UNIVERSITY, SEOUL 02841, S. KOREA
*Email address*: `yaryong@korea.ac.kr`

SUNGGEUM HONG, DEPARTMENT OF MATHEMATICS, CHOSUN UNIVERSITY, GWANGJU 61452, S. KOREA
*Email address*: `skhong@chosun.ac.kr`

CHAN WOO YANG, DEPARTMENT OF MATHEMATICS, KOREA UNIVERSITY, SEOUL 02841, S. KOREA
*Email address*: `cw_yang@korea.ac.kr`